\documentclass[a4paper,12pt]{article}

\usepackage{amsmath,amsfonts,latexsym,theorem,%showkeys
}
\usepackage{graphicx}
\usepackage{epsfig}
\usepackage{epstopdf}
\usepackage[ruled,vlined]{algorithm2e}
\numberwithin{equation}{section}

\newtheorem{Theorem}{Theorem}[section]
\newtheorem{Definition}{Definition}[section]
\newtheorem{Proposition}{Proposition}[section]
\newtheorem{Lemma}{Lemma}[section]

\newtheorem{Corollary}{Corollary}[section]
\newenvironment{Proofc}[1]{\smallskip\par\noindent\textsc{#1}\quad}%
  {\hfill$\Box$\bigskip\par}
\newenvironment{Proof}{\begin{Proofc}{Proof}}{\end{Proofc}}

\newtheorem{Remark}{Remark}[section]
\def\a{\alpha}

\def\D{\Delta}
\def\g{\gamma}
\def\G{\Gamma}
\def\l{\lambda}

\def\s{\sigma}

\def\th{\theta}

\def\t{\tau}
\def\e{\varepsilon}

\def\pd{\partial}
\def\ds{\displaystyle}
\def\half{\frac{1}{2}}
\newcommand{\R}{{\mathbb R}}
\newcommand{\N}{{\mathbb N}}

\def\pd{\partial}

\begin{document}
\title{A  model problem for Mean Field Games on networks\footnotemark[1]}
\date{}

%\author{Fabio Camilli , Elisabetta Carlini, Claudio Marchi}
%\subjclass{Primary:  ; Secondary: .}
%\keywords{.}

%\email{camilli@dmmm.uniroma1.it}
%\email{}
%\email{}
%%%%%%%%%%%%%%%%%%%%%%%%%%%%%%%%%%%%%%%%%%%%%%%%%%%%%%%%%%%%%%%%
\maketitle
\vskip -1.5cm
\centerline{\scshape Fabio Camilli}
\medskip
{\footnotesize
\centerline{``Sapienza" Universit{\`a}  di Roma}
\centerline{Dip. di Scienze di Base e Applicate per l'Ingegneria}
\centerline{via Scarpa 16, 0161 Roma, Italy}
\centerline{{\tt e-mail:camilli@dmmm.uniroma1.it}}
}
\medskip
\centerline{\scshape Elisabetta Carlini}
\medskip
{\footnotesize
\centerline{``Sapienza" Universit{\`a}  di Roma}
\centerline{Dip. di Matematica}
\centerline{P.le A. Moro 5, 00185 Roma, Italy}
\centerline{{\tt e-mail:carlini@mat.uniroma1.it}}
}
\medskip
\centerline{\scshape Claudio Marchi}
\medskip
{\footnotesize
\centerline{Universit\`a di Padova}
\centerline{Dip. di Matematica }
\centerline{via Trieste 63, 35121 Padova, Italy}
\centerline{{\tt email:marchi@math.unipd.it}}
}
\bigskip
\footnotetext[1]{This work is partially supported by European Union under the 7th Framework Programme FP7-PEOPLE-
2010-ITN Grant agreement number 264735-SADCO}

\begin{abstract}
In \cite{gll}, Gu{\'e}ant, Lasry and Lions considered the  model problem ``What time does meeting start?'' as a prototype for a general class of optimization problems with a continuum of players, called Mean Field Games problems. In this paper we  consider a similar model, but with the dynamics of  the agents   defined   on a network. We discuss appropriate transition conditions at the vertices which give  a well posed    problem   and we   present some numerical results.
\end{abstract}
 \begin{description}
\item [{\bf MSC 2000}:] 91A15, 35R02, 35B30, 49N70, 65M06. %93N20.
 \item [{\bf Keywords}:] networks,  mean field games, stochastic optimal control,  numerical methods.
 \end{description}
%%%%%%%%%%%%%%%%%%%%%%%%%%%%%%%%%
%%%        Introduzione        %%
%%%%%%%%%%%%%%%%%%%%%%%%%%%%%%%%%
%
\section{Introduction}
The study of pedestrian flow in a crowd environment is attracting an increasing interest and    some models based on   optimization principles have been  recently proposed, see for example   \cite{bdmw,do,pt}. In  some applications (crowd motion  in shopping centers, stations, airports)   the dynamics of the population is defined on   a network  rather than in an Euclidean domain.\\
There is a large literature concerning   vehicular traffic on   road networks (see \cite{gp} and reference therein). These models are   based on a fluido-dynamical approach with the dynamics  described by some nonlinear conservation law and  appropriate transition conditions at the junctions modelling the interactions of the cars coming from different roads.
Vehicular traffic models do not seem to
be  adequate  to reproduce the pedestrian flow since they do not take into account  the interactions and the goal-directed decisions of the agents.\\
 Aim of this paper is to study a simple optimization model for  the evolution of a large number  of agents moving on a network. The model is based on the one described in \cite{gll},  titled ``What time does meeting start?'', and consists in finding the optimal arrival time at a place where the meeting is being held with the  starting time  defined by means of  a quorum rule.
This problem  can be considered as a prototype for a large class of optimization problems based on the Mean Field Game  (MFG) theory.
This theory has been   introduced by Lasry and Lions \cite{ll} (see also \cite{a}, \cite{c}, \cite{gs}) with the  aim of describing the behavior of very large number of agents who take decisions in a context of strategic interactions.
%Let us recall that in \cite{a} Achdou studied (also from a numerical point of
%view) the case where a population is driven to leave a closed room with
%obstacles. Our problem is different from its one because our equations are
%coupled via the boundary conditions and, mainly, they are defined on a
%network.
\\
The main difficulties in our approach is  to deal with the transition conditions at the  internal vertices to obtain a well posed MFG problem.
It is  known that a parabolic equation on a network has to be complemented  with the usual initial-boundary conditions and some transition  conditions at the internal vertices (see \cite{vbn,mu}).
In fact, in  our model the stochastic differential equation describing the motion of the agent  inside the arcs is coupled with a condition  prescribing the probability that it  enter in a given  edge when it occupies a transition vertex; this fact give rise to a Kirchhoff type condition (see \cite{fs}).
 Using an appropriate change of variable we transform the original MFG system in a forward-backward system of two heat equations
coupled via the initial datum. Relying on classical results for the heat equation on   networks and some appropriate estimates
for the specific problem, we prove the well-posedness of the heat system and the existence of a mean field for the quorum problem.\par
Going back to the original MFG problem  we obtain existence and uniqueness of the solution to a system composed by a {\it backward} Hamilton-Jacobi-Bellman equation and a {\it forward} Fokker-Planck on the arcs with  transitions conditions  expressing respectively the probability that a single agent enters  a given arc    and  the conservation of the density of the agents  through a vertex.\par
The paper is organized as follows. In Section \ref{S2} we describe the model problem. In \ref{S3} we prove some technical results
concerning the heat equation on the network which are used in Section \ref{S4} to show the existence of the mean field. in Section \ref{S5}, we illustrate the problem with some numerical examples. Finally, the Appendix contains some technical proofs.

\vskip 8pt
\par\textbf{Notations:}
A  network   is a finite  collection of points $V:=\{v_i\}_{i\in I}$ in $\R^n$ connected by continuous, non self-intersecting arcs $E:=\{e_j\}_{j\in J}$. Each arc $e_j$ is  parametrized by a smooth function $\pi_j:[0,l_j]\to\R^n,\, l_j>0$.\\
For $i\in I$ we set $Inc_i:=\{j\in J\mid\,e_j \,\text{is incident to}\,v_i\}$. We denote by  $I_B:=\{i\in I\mid \# Inc_i=1\}$, $I_T:=I\setminus I_B$, by  $\partial \G:=\{v_i\in V\mid \,i\in I_B\}$, the set of boundary vertices of $\G$, and by $\G_T:=\{v_i\mid i\in I_T\}$, the set of transition vertices.\\
The network is not oriented, but the parametrization of the arcs induces an orientation which can be expressed by the \emph{signed incidence matrix} $A=\{a_{ij}\}$ with
\begin{equation*}\label{incidence}
   a_{ij}:=\left\{
            \begin{array}{rl}
              1 & \hbox{if $v_i\in  e_j$ and $\pi_j(0)=v_i$,} \\
              -1 & \hbox{if $v_i\in  e_j$ and $\pi_j(l_j)=v_i$,} \\
              0 & \hbox{otherwise.}
            \end{array}
          \right.
\end{equation*}
In the following we always identify   $x\in   e_j$ with   $y=\pi_j^{-1}(x)\in [0,l_j]$. For any function $u:  \G\to\R$ and each $j\in J$ we denote by $u_j:[0,l_j]\to \R$ the restriction of $u$ to $e_j$, i.e.  $u_j(y)=u(\pi_j(y))$ for $y\in [0,l_j]$.
For $\g \in\N$, we define differentiation along an edge $e_j$ by
\[\pd^\g_{j}u(x):= \frac{d^\g u^j}{d y^\g}  (y),\qquad\text{for $y=\pi^{-1}_j(x)$, $x\in e_j$}\]
and at a vertex $v_i$ by
\[\pd^\g_{j}u(v_i):= \frac{d^\g u^j}{d y^\g}  (y)\qquad\text{for $y=\pi^{-1}_j(v_i)$,  $j\in Inc_i$.}\]

%%%%%%%%%%%
%         %
%%%%%%%%%%%
\section{The model problem}\label{S2}
Following \cite{gll}, we   describe  the   model ``What time does meeting start?'' with the variant that the dynamics of the agents are defined on a network $\G$.
For the sake of simplicity, we assume that the place where the meeting  is being held is the unique boundary vertex, namely $\pd \G=\{v_0\}$; the general case can be dealt with by using easy adaptations.
The meeting is scheduled at a certain time $t_0$ but the common experience says that in general it starts at a time $T$ {\it greater than} $t_0$, when a certain rule is reached, for example the presence of a certain percentage of participants.\\
At the initial time there is a
 continuum of indistinguishable players distributed according to a distribution function $m_0:\G\to\R$.
The player's dynamics is subject to random perturbations. We assume that, inside each edge $e_j$, the generic  agent moves according to the process
 \begin{equation}\label{sdetoy}
    dX (t)=a (t)dt+ \s dW (t)
 \end{equation}
where the drift $a$ is the  control  variable (and it coincides with the speed), $\s =(\s_j)_{j\in J}$ with $\s_j>0$ and  $W$ is a Brownian process, which is an  independent disturbance  for each player. Moreover we assume that, at each transition vertex $v_i$, it spends zero time a.s. and it enters in one of the  incident edge~$e_j$ with    probability $ 1/ \#(Inc_i)$  (see \cite{fs,fw} for stochastic differential equations on networks). We  denote by $\t$ the  random time the agent  reaches $v_0$, i.e.
\[
\tau:=\inf\{t>0:\,X(t)\in\pd\G\}.\]
Moreover each player wants to optimize its arrival time $\tau$ taking into account various parameters, which are encoded
in the cost functional
\begin{equation}\label{costtoy}
 J(x,t, a(\cdot))=   \int_t^{\tau\wedge T_{max}}  \half a^2(t)dt+c(\t\wedge T_{max})
\end{equation}
where $\half a^2(t)$ is the actual cost of moving along the network at the velocity $a$ while $c$ is the final cost and $T_{max}\in \R$ is
a time which cannot be exceeded for the end of the meeting. The cost function $c:[0, T_{max}]\to \R$ is given by
\begin{equation}\label{boundarycost}
c(s)=c_1(s-t_0)+c_2(s-T)+c_3(T-s),\qquad s\in [0,T_{max}]
\end{equation}
where    $c_i:\R\to\R$, $i=1,2,3$ are smooth functions  such that $c_i(s)=0$ for $s\le 0$ and $c_i(s)>0$ for $s>0$.
 The term $c_1(s-t_0)$ represents a reputation cost of lateness in relation to scheduled time $t_0$;
the term $c_2(s-T)$ a cost of lateness in relation to actual starting time of the meeting $T$; $c_3(T-s)$ a waiting time cost which corresponds to the time lost waiting the starting of the meeting.
It is worth noticing   that the cost $c$ depends on $T$  via the cost of lateness and the cost of waiting; hence, in order to display this dependence, from now on  we write $c_T$.\par
Nash equilibrium theory assumes that each player want to optimize the arrival time by assuming that
actual time  $T$ the meeting starts  is known. Hence  each agent has to solve the optimization
problem
\begin{equation}\label{valuetoy}
  u(x,t) = \min_{a(\cdot)}  J(x,t, a(\cdot))
\end{equation}
where $(x,t)\in \G\times [0,T_{max}]$.
Note that $\max_{a\in \R}\{-ap+\half|a|^2\}=-\half |p|^2$ for $a=-p$ and the optimal control in feed-back form is given by $a^*(x,t)=-\pd_xu(x,t)$.
By an application of the Dynamic Programming Principle the value function,  if it is assumed to be regular, formally solves  the  Hamilton-Jacobi-Bellman equation
\begin{equation*}%\label{hjbtoy1}
    \ds  \pd_t u+\nu\pd^2_{x}u+\half|\pd_x u|^2=0\quad (x,t)\in \G\times (0,T_{max}),
\end{equation*}
where $\nu=\s^2/2$ (i.e., $\nu_j=\s_j^2/2 \quad \forall j\in J$),
with final-boundary conditions and transition on internal vertices (Kirchhoff condition)
\begin{align*}
&u(x,T_{max})=c_T(T_{max})\quad x\in\G, \qquad
u(v_0,t)=c_T(t)\quad s\in[0, T_{max}],\\
&\sum_{j\in Inc_i} a_{ij} \pd_{j} u(v_i ,s)=0\qquad (v_i,s)\in \G_T\times (0,T_{max}).
\end{align*}
%\begin{align}
%    u(v_0,t)&=c_T(t)\quad &s\in[0, T_{max}],\label{hjbtoy2}\\
%   u(x,T_{max})&=c_T(T_{max})&x\in\G.\label{hjbtoy3}
%\end{align}

On the other hand, by duality, the dynamic of the agents, i.e. the evolution of the initial distribution $m_0$, is governed inside each edge by the Fokker-Planck equation
\begin{equation*}%\label{fptoy}
    \pd_t m-\nu\pd_x^2m-\pd_x((-\pd_x u)m)=0\quad (x,t)\in \G\times (0,T_{max})
\end{equation*}
and we assume the initial-boundary condition (with a ``smooth fit'') and a Kirchhoff condition on internal vertices
\begin{align*}
&  m(x,0)=m_0(x)\quad x\in\G, \qquad m(v_0,s)=0 \qquad  s\in [0,T]\\
&\sum_{j\in Inc_i} a_{ij}\nu_{j}[ \pd_{j}m - m\pd_{j}u ](v_i,s)=0\qquad(v_i,s)\in \G_T\times (0,T_{max}).
\end{align*}
Observe that the previous  Kirchhoff condition  implies that the parabolic  flux of the agents is null at the junctions, giving the conservation of the total mass (see \cite{cg} for similar assumptions).
%It is worth noticing that the function $u$ and consequently the function $m$ via   the drift term $-\pd_x u$   depend on the actual time $T$.\\

The flow of participants reaching $v_0$ is given by $s\mapsto \pd_x m(v_0,s)$, hence  the cumulative distribution $F$ of the arrival times
is \[F(s)=\int_0^s \nu \pd_x m(v_0,r)dr.\]
The actual starting time  $T$ is fixed by a quorum rule, which means that the meeting starts when a given percentage $\theta$  of the participants  has reached the meeting place $v_0$. Given $m$, we set
\begin{equation}\label{quorum}
   T=\left\{
       \begin{array}{ll}
         t_0, & \hbox{$F^{-1}(\th)\le t_0 $ } \\[4pt]
         F^{-1}(\th), & \hbox{$t_0<F^{-1}(\th)< T_{max}$}\\[4pt]
         T_{max}, & \hbox{$F^{-1}(\th)\ge T_{max}$.}
       \end{array}
     \right.
\end{equation}
Note that $T$ is the mean field, i.e. the information that the single agent  has about the behavior of the other agents:
the starting rule induces a strategic interactions among the participants and  $T$ influences  as
an external field the decisions of the   agents. The main point is to prove the existence and the uniqueness of a time $T$ which is coherent with the expectations of the participants. As in \cite{gll}, this can be done by proving that the scheme:
\begin{equation}\label{scheme}
T\rightarrow u\rightarrow m\rightarrow T^*
\end{equation}
with $T^*$  defined by \eqref{quorum}, has a fixed point in $[t_0,T_{max}]$. To this end, it is important to study  existence and uniqueness of a solution to
the  forward-backward system
\begin{equation}\label{MFGtoy}
\left\{\begin{array}{ll}
\ds \pd_t u+\nu\pd^2_xu+\half|\pd_x u|^2=0\quad &(x,s)\in \G\times (0,T_{max})\\[6pt]
\pd_t m-\nu\pd^2_x m+\pd_x( \pd_x u \,m)=0  &(x,s)\in \G\times (0,T_{max})\\[6pt]
\sum_{j\in Inc_i} a_{ij}  \pd_{j} u(v_i ,s)=0\quad &(v_i,s)\in \G_T\times (0,T_{max})\\[6pt]
\sum_{j\in Inc_i} a_{ij}\nu_{j}[\pd_{j}m - m\pd_{j}u ](v_i,s)=0 \quad &(v_i,s)\in \G_T\times (0,T_{max})\\[6pt]
m(x,0)=m_0(x),\, u(x,T_{max})=c_T(T_{max})&x\in\G\\[6pt]
m(v_0,s)=0,\, u(v_0,s)=c_T(s)   & s\in[0, T_{max}].
\end{array}\right.
\end{equation}
For the sake of simplicity, from now on we assume
\begin{equation}\label{constants}
\nu_j=1\qquad \forall j\in Inc_i.
\end{equation}
As in \cite{gll, g}  we apply a change of variable which transforms system~\eqref{MFGtoy} into a forward-backward system of heat equations coupled through the initial conditions.

\begin{Proposition} If $(\phi,\psi)$ is  a smooth solution of the system
\begin{equation}\label{MFGchange}
  \left\{
    \begin{array}{ll}
     -\pd_t \phi-\pd^2_x\phi=0\quad &(x,s)\in \G\times (0,T_{max}),\\[6pt]
     \pd_t \psi-\pd^2_x\psi=0       &(x,s)\in \G\times (0,T_{max})\\[6pt]
     \sum_{j\in Inc_i} a_{ij} \pd_{j} \phi(v_i ,s)=0& (v_i,s)\in \G_T\times (0,T_{max})\\[6pt]
     \sum_{j\in Inc_i} a_{ij} \pd_{j} \psi(v_i ,s)=0& (v_i,s)\in \G_T\times (0,T_{max})\\[6pt]
     \psi(x,0)=\frac{m_0(x)}{\phi(x,0)},\, \phi(x,T_{max})=e^{c_T(T_{max})}   &x \in\G\\[6pt]
     \psi(v_0,s)=0,\, \phi(v_0,s)=e^{c_T(s)}& s\in[0, T_{max}]
   \end{array}
 \right.
\end{equation}
with $\phi>0$, then
\begin{equation}\label{cov}
(u,m)=(\ln(\phi), \phi\,\psi)
\end{equation}
is a solution of system~\eqref{MFGtoy}.
\end{Proposition}
\begin{Proof}
Let $(\phi,\psi)$ and $(u,m)$ be defined as in the statement.
The proofs that $(u,m)$ is a solution to the PDEs and to initial-final-boundary conditions of \eqref{MFGtoy} follow by easy calculations; hence, we shall omit them.
Let us prove that $(u,m)$ verifies the transitions condition of \eqref{MFGtoy}. Since $\phi=e^u$, we get
\[ 0=\sum_{j\in Inc_i} a_{ij} \pd_{j} \phi=e^u\sum_{j\in Inc_i} a_{ij} \pd_{j} u\]
which amounts to the first transition condition in~\eqref{MFGtoy}. On the other hand, since $\psi=m e^{-u}$, we have
\[ 0=\sum_{j\in Inc_i} a_{ij} \pd_{j} \psi=e^{-u}\sum_{j\in Inc_i} a_{ij} (\pd_{j} m-m\pd_{j} u).\]
Taking into account the previous relation, we obtain the second transition condition in~\eqref{MFGtoy}.
\end{Proof}
\begin{Remark}
It is worth to observe that, by similar arguments, one can linearize a more general class of MFG systems (see \cite{g}). Actually, assume that $\nu_j$   are  positive constants and that the cost $J$ in \eqref{costtoy} includes a potential term depending on the distribution of other players, i.e.
\[
 J(x,t, a(\cdot))=   \int_t^{\tau\wedge T_{max}} \big[ \half a^2(t)+f(X(t),m(t))\big]dt+c(\t\wedge T_{max}).
\]
In this case, in the system~\eqref{MFGtoy} the Hamilton-Jacobi-Bellman equation is
\begin{equation*}
\pd_t u+\nu\pd^2_xu+\half|\pd_x u|^2=-f(x,m)\quad (x,s)\in \G\times (0,T_{max}),
\end{equation*}
while the Fokker-Planck equation   and the boundary-transition conditions are left unchanged. Now, $(\phi,\psi)=(e^{u/\s^2}, me^{-u/\s^2})$ solve
\begin{equation*}
\left\{
\begin{array}{l}
-\pd_t \phi-\nu\pd^2_x\phi=-\frac\phi{2\nu} f(x,\phi\psi),\quad \pd_t \psi-\nu\pd^2_x\psi=\frac\psi{2\nu} f(x,\phi\psi)\quad \textrm{in }\G\times (0,T_{max}),\\[7pt]
\sum\limits_{j\in Inc_i} a_{ij}  \pd_{j} \phi(v_i ,s)= \sum\limits_{j\in Inc_i} a_{ij}\nu_j(\phi\pd_{j} \psi)(v_i ,s)=0 \quad \textrm{in }\G_T\times (0,T_{max})\\[7pt]
\psi(\cdot,0)=\frac{m_0(\cdot)}{\phi(\cdot,0)},\, \phi(\cdot,T_{max})=e^{\frac{c_T(T_{max})}{\s^2}},\, \psi(v_0,\cdot)=0,\, \phi(v_0,\cdot)=e^{\frac{c_T(\cdot)}{\s^2}}.
   \end{array} \right.
\end{equation*}
\end{Remark}
\section{The heat equation on a network}\label{S3}
In this section we collect some technical results about existence, uniqueness and  a priori estimates for classical solutions to \eqref{MFGchange}. These results will be used in the next section to prove the existence of the mean field $T$. \\
%Note that the system is composed of two heat equations coupled only in the initial condition for $\psi$.\\
We introduce some functional spaces on the network.
%A function  $u$ is  continuous on $\G$ and we write $u\in C(\G)$  if it is continuous with respect to the subspace topology of $\G$, namely, $u^j\in C([0,l_j])$ for any $j\in J$ and $u^j(\pi_j^{-1}(v_i))=u^k(\pi_k^{-1}(v_i))$  for any $i\in I$, $j,k\in Inc_i$. \par
%For $q\in\N$, we set
%\begin{equation*}
%C^{q}(\G):=\{u\in C(\G)\mid \forall j\in J,\,u_j\in C^{q}([0,l_j])\}
%\end{equation*}
%which is a Banach space with respect to its norm \[|u|^{(q)}_{\G}:=\sup_{\g\le q}\big[\sup_{j\in J}|\pd_j^\g u|_{L^\infty(0,l_j)}\big]\]
%(observe that no continuity condition at the     vertices  is  prescribed for the derivates   of a  function $u\in C^k(\G)$).
We recall that the Sobolev space $W^{2,1}_{q,(a,b)\times(0,T)}$ (with $q\geq 1$) consists  of the elements of $L^q((a,b)\times(0,T))$ having generalized derivatives of the form $\pd_t^r\pd_x^s$ with $2r+s=2$ and it is endowed with its usual norm (see \cite{lsu}). %; its norm is defined as
%\[|u|^{2,1}_{q,(a,b)\times(0,T)}:=\sum_{j=0}^2 \sum_{2r+s=j}|\pd_t^r\pd_x^su|_{L^q((a,b)\times(0,T))}.\]
For $q\in\N$ and $\a\in(0,1)$, $C^{(q+\a)}([a,b])$ stands for the Banach space of $q$ times differentiable functions on $[a,b]$, whose $q$-th derivative is H\"older continuous with exponent $\a$ and it is endowed with the usual H\"older norm $|\cdot|^{(q+\a)}_{[a,b]}$. For $\a\in(0,1)$,  $C^{(2+\a,1+\a/2)}([a,b]\times[0,T])$, with the norm $|\cdot|^{(2+\a,1+\a/2)}_{[a,b]\times[0,T]}$, denotes the Banach space of functions $f:[a,b]\times[0,T]\to \R$ which have  Holder continuous derivatives $\pd^2_x f$ and $\pd_tf$.
%and finite norm
%\begin{multline*}
%|f|^{(2+\a,1+\a/2)}_{[a,b]\times[0,T]}:=
%|f|_\infty+|f_t|_\infty+|f_x|_\infty+|f_{xx}|_\infty+
%H^\a_{x}(f_{xx})+H^{\a/2}_{t}(f_{t})+H^{(1+\a)/2}_{x}(f_{t})\\+
%H^{\a/2}_{t}(f_{xx})
%\end{multline*}
%where the last four summands denote the H\"older constants with respect to the indicate exponents and variables.
%%
\begin{Definition}\hfill
\begin{itemize}
\item[i)] For $q\in\N$ and $\a\in(0,1)$, we set
\begin{equation*}
C^{(q+\a)}(\G):=\{u\in C(\G)\mid \forall j\in J,\,u_j\in C^{(q+\a)}([0,l_j])\}
\end{equation*}
which is a Banach space with respect to its norm $|u|^{(q+\a)}_{\G}:=\sup_{j\in J}|u_j|^{(q+\a)}_{[0,l_j]}$.
\item[ii)]  For $\a\in(0,1)$, we set
\begin{equation*}
C^{(2+\a,1+\a/2)}(\G\times[0,T]):=\{u\in C(\G\times[0,T])\mid \forall j\in J,\,u_j\in C^{(2+\a,1+\a/2)}([0,l_j]\times[0,T])\}
\end{equation*}
which  is a Banach space with respect to its norm $|u|^{(2+\a,1+\a/2)}(\G\times[0,T]):=\sup_{j\in J}|u_j|^{(2+\a,1+\a/2)}_{[0,l_j]\times[0,T]}$.
\end{itemize}
\end{Definition}
%%

%We shall also need some Banach spaces of H\"older continuous function on the graph.For $q\in\N$ and $\a\in(0,1)$, $C^{(q+\a)}([a,b])$ stands for the Banach space of $q$ times differentiable functions on $[a,b]$, whose $q$-th derivative is H\"older continuous with exponent $\a$ and it is endowed with the usual H\"older norm $|\cdot|^{(q+\a)}_{[a,b]}$. For $q\in\N$ and $\a\in(0,1)$, we set
%\begin{equation*}
%C^{(q+\a)}(\G):=\{u\in C(\G)\mid \forall j\in J,\,u_j\in C^{(q+\a)}([0,l_j])\}
%\end{equation*}
%which is a Banach space with respect to its norm $|u|^{(q+\a)}_{\G}:=\sum_{j\in J}|u_j|^{(q+\a)}_{[0,l_j]}$.
%
%Moreover, for $\a\in(0,1)$, we introduce the space
%\begin{equation*}
%C^{(2+\a,1+\a/2)}(\G\times[0,T]):=\{u\in C(\G\times[0,T])\mid \forall j\in J,\,u_j\in C^{(2+\a,1+\a/2)}([0,l_j]\times[0,T])\}
%\end{equation*}
%where $C^{(2+\a,1+\a/2)}([a,b]\times[0,T])$, with the norm $|\cdot|^{(2+\a,1+\a/2)}_{[a,b]\times[0,T]}$, denotes the Banach space of functions $f:[a,b]\times[0,T]\to \R$ which have continuous derivatives $f_{xx}$ and $f_t$ and finite norm
%\begin{multline*}
%|f|^{(2+\a,1+\a/2)}_{[a,b]\times[0,T]}:=
%|f|_\infty+|f_t|_\infty+|f_x|_\infty+|f_{xx}|_\infty+
%H^\a_{x}(f_{xx})+H^{\a/2}_{t}(f_{t})+H^{(1+\a)/2}_{x}(f_{t})\\+
%H^{\a/2}_{t}(f_{xx})
%\end{multline*}
%where the last four summands denote the H\"older contants with respect to the indicate exponents and variables. Finally, the space $C^{(2+\a,1+\a/2)}(\G\times[0,T])$ is a Banach space when endowed with the norm $|u|^{(2+\a,1+\a/2)}(\G\times[0,T]):=\sum_{j\in J}|u_j|^{(2+\a,1+\a/2)}_{[a,b]\times[0,T]}$.

In the next proposition we establish the well-posedness of the initial-boundary problem for the heat equation obtained by the Hamilton-Jacobi-Bellman equation of \eqref{MFGtoy} via the change of variable \eqref{cov}.
\begin{Proposition}\label{exisheat}
Assume that $w_0\in C^{(1+\a/2)}([0,T_{max}])$, for some $\a\in(0,1)$.   Then there exists a unique solution $w\in C^{(2+\a,1+\a/2)}(\G\times[0,T_{max}])$ of the problem
\begin{equation}\label{heat}
\ds\left\{
\begin{array}{ll}
-\pd_t w- \pd^2_xw=0\quad& (x,s)\in \G\times (0,T_{max})\\[6pt]
\sum_{j\in Inc_i} a_{ij}\pd_{j} w(v_i ,s)=0& (v_i,s)\in \G_T\times (0,T_{max})\\[6pt]
w(v_0,s)=w_0(s)& s\in[0, T_{max}] \\[6pt]
w(x,T_{max})=  w_0(T_{max}) & x\in\G.
\end{array}\right.
\end{equation}
%where $A=\{a_{ij}\}_{ij}$ is  the incidence matrix defined in \eqref{incidence}.
Moreover, the following estimate holds
\begin{equation}\label{heatvB}
|w|^{(2+\a,1+\a/2)}_{\G\times[0,T_{max}]}\leq K_0 |w_0|^{(1+\a/2)}_{[0,T_{max}]}
\end{equation}
where $K_0$ is a constant independent of $w_0$.
Finally, for $w_0>0$, we have $w\geq \min w_0$ in  $\G\times[0,T_{max}]$.
\end{Proposition}
\begin{Proof}
The statement is an immediate consequence of the result in \cite{vb}. Let us just note that the compatibility conditions in \cite{vb}  are obviously satisfied because the terminal condition is constant and the right-hand side of the Kirchhoff condition is null.
Moreover the strict positivity of $w$ is a consequence of the comparison principle for classical solution of the heat equation (see \cite{vbn}). We observe that it can be proved using the same arguments of \cite[Proposition 2]{g}.
\end{Proof}
Since $v_0$ is a boundary vertex,   there exists a unique edge, say $e_0$  incident to it. Without any loss of generality, we denote $v_1$ the other endpoint of $e_0$ and we assume that the parametrization of $e_0$ fulfills:
\begin{equation}\label{v1}
\pi_0(0)=v_0\quad \textrm{and } \pi(l_0)=v_1.
\end{equation}
For $\l\in(0,1)$, we set
\begin{equation}\label{vprimolambda}
e_{0,\l}:=\pi_0([0,\l l_0]), \qquad v'_\l:=\pi_0(\l l_0)
\end{equation}
namely, $v'_\l$ is a point in the edge $e_0$ while $e_{0,\l}$ is the part of $e_0$ between $v_0$ and $v'_\l$.

In the next proposition, we establish existence and uniqueness of a classical solution to the heat equation obtained by the Fokker-Planck equation of
\eqref{MFGtoy} via \eqref{cov}. Moreover we show a ``weak" continuous dependence estimate  in the sub-edge $e_{0,1/2}$ with respect to the initial datum $\mu(\cdot)/w(\cdot,0)$ where $w$ is the solution of \eqref{heat}.

\begin{Proposition}\label{exisFP0}% (Existence for the FP eq - 1)
Let  $w$ be the solution of  problem \eqref{heat} and assume
\begin{equation}\label{33weak1}
\mu_0\in C^{(2+\a)}(\G),\quad\textrm{with }
\mu_0(v_0)=0.
\end{equation}
 Then there exists a unique solution $\mu\in C^{2,1}(\G\times(0,T_{max}))\cap C^0(\bar \G\times[0,T_{max}])$ of the problem
\begin{equation}\label{FPnet}
\ds\left\{
\begin{array}{ll}
\pd_t \mu- \pd^2_x\mu=0\quad& (x,s)\in \G\times (0,T_{max})\\[6pt]
\sum_{j\in Inc_i} a_{ij} \pd_{j} \mu(v_i ,s)=0& (v_i,s)\in \G_T\times (0,T_{max})\\[6pt]
\mu(v_0,s)=0& s\in[0, T_{max}] \\[6pt]
\mu(x,0)= \frac{\mu_0(x)}{w(x,0)} & x\in\G.
\end{array}\right.
\end{equation}
Moreover, for every $q\geq 1$, the following estimate holds
\begin{equation}\label{33weak}
|\mu|^{2,1}_{q,e_{0,1/2}\times[0,T_{max}]}\leq K_1 |\mu_0/w(\cdot,0)|^{(2+\a)}_{\G}
\end{equation}
where $K_1$ is a constant independent of $\mu_0$ and $w$.
\end{Proposition}
The proof is postponed in the Appendix.
%%%%%%%%

In the next proposition, we establish two continuous dependence estimates for the solution of problem \eqref{FPnet} with respect to the initial datum: the former is a ``strong'' estimate in the sub-edge $e_{0,1/2}$ while the latter is the classical estimate in the whole network.
%
%In the next Proposition, we establish a ``strong'' continuous dependence estimate in the sub-edge $e_{0,1/2}$ and a classical continuous dependence estimate in the whole network with respect to the initial datum for the solution of problem \eqref{FPnet}.
\begin{Proposition}\label{exisFP} %(Existence for the FP eq - 2)
Let $\mu$ be the solution to \eqref{FPnet}. Besides the hypotheses of  Proposition~\ref{exisFP0}, assume
\begin{equation}\label{cmpv0}
\pd_x \mu_0(v_0)=\pd^2_x \mu_0(v_0)=0.
\end{equation}
% Then there exists a unique solution $\mu\in C^{2,1}(\G\times(0,T_{max}))\cap C^0(\bar \G\times[0,T_{max}])$ of the problem
%\begin{equation}\label{FPnet}
%\ds\left\{
%\begin{array}{ll}
%\pd_t \mu- \nu\,\pd^2_x\mu=0\quad& (x,s)\in \G\times (0,T_{max})\\[6pt]
%\sum_{j\in Inc_v} a_{ij}\rho_{ij} \pd_{x,j} \mu(v_i ,s)=0& (v_i,s)\in \G_T\times (0,T_{max})\\[6pt]
%\mu(v_0,s)=0& s\in[0, T_{max}] \\[6pt]
%\mu(x,0)= \frac{\mu_0(x)}{w(x,0)} & x\in\G
%\end{array}\right.
%\end{equation}
%where $A=\{a_{ij}\}_{ij}$ is  the incidence matrix defined in \eqref{incidence}.
\begin{itemize}
\item[i)] There holds
\begin{equation}\label{IV10.1}
|\mu|^{(2+\a,1+\a/2)}_{e_{0,1/2}\times[0,T_{max}]}\leq K_2 |\mu_0/w(\cdot,0)|^{(2+\a)}_{\G}
\end{equation}
where $K_2$ is a constant independent of $\mu_0$ and $w$.
\item[ii)] Under the further assumption
\begin{equation}\label{FPallaVB}
\pd_j \mu_0(v_i)=\pd^2_j\mu_0(v_i)=0 \qquad \forall i\in I_T,\, j\in Inc_i,
\end{equation}
the function $\mu$ belongs to $C^{(2+\a,1+\a/2)}(\bar \G\times[0,T_{max}])$ and verifies
\[
|\mu|^{(2+\a,1+\a/2)}_{\G\times[0,T_{max}]}\leq K_3 |\mu_0/w(\cdot,0)|^{(2+\a)}_{\G}
\]
where $K_3$ is a constant independent of $\mu_0$ and $w$.
\end{itemize}
\end{Proposition}
%In the next Proposition, we shall obtain the regular continuous dependence in the whole network $\G$ of the solution to the Fokker-Planck problem; to this end we shall need some compatibility conditions for the initial datum $\mu_0$ stronger than those in \eqref{cmpv0}.
%
%\begin{Proposition}\label{exisFP2} (Existence for the FP eq -- 3)
%Assume the hypotheses of Proposition \ref{exisFP0}, relation \eqref{cmpv0} and

%\end{Proposition}
%\begin{Proof}
%
%\end{Proof}
The proof is postponed in the Appendix. Let us now establish a well-posedness result for the system \eqref{MFGchange}.

\begin{Theorem}\label{thm:MFG}
Assume that, for some $\a\in(0,1)$, there holds
\begin{equation}\label{ipminime}
c_T \in C^{(1+\a/2)}([0,T_{max}]),\quad c\geq 0,\qquad m_0\in C^{(2+\a)}(\G)\textrm{ with }  m_0(v_0)=0.
\end{equation}
Then, there exists a unique classical solution $(\phi,\psi)$ to the  system~\eqref{MFGchange}
%\begin{equation}\label{MFGchangefull}
%  \left\{
%    \begin{array}{ll}
%     -\pd_t \phi-\nu\pd^2_x\phi=0\quad &(x,s)\in \G\times (0,T_{max}),\\[6pt]
%     \pd_t \psi-\nu\pd^2_x\psi=0       &(x,s)\in \G\times (0,T_{max})\\[6pt]
%     \sum_{j\in Inc_v} a_{ij}\nu_j\rho_{ij} \pd_{x,j} \phi(v_i ,s)=0& (v_i,s)\in \G_T\times (0,T_{max})\\[6pt]
%     \sum_{j\in Inc_v} a_{ij}\nu_j\rho_{ij} \pd_{x,j} \psi(v_i ,s)=0& (v_i,s)\in \G_T\times (0,T_{max})\\[6pt]
%     \psi(x,0)=\frac{m_0(x)}{\phi(x,0)},\, \phi(x,T_{max})=e^{ \frac{c_T(T_{max})}{\s^2}}   &x \in\G\\[6pt]
%     \psi(v_0,s)=0,\, \phi(v_0,s)=e^{ \frac{c_T(s)}{\s^2}}& s\in[0, T_{max}]
%   \end{array}
% \right.
%\end{equation}
with $\phi>0$.
Moreover, the following estimates hold
\begin{equation*}
\begin{array}{rl}
(i)& \phi\geq 1, \quad
|\phi|^{(2+\a,1+\a/2)}_{\G\times[0,T_{max}]}\leq K |c_T|^{(1+\a/2)}_{[0,T_{max}]},\quad
|\psi|^{2,1}_{q,e_{0,1/2}\times[0,T_{max}]}\leq K |m_0/\phi(\cdot,0)|^{(2+\a)}_{\G}\\[12pt]
(ii)&\textrm{If $m_0$ fulfills \eqref{cmpv0}:}\qquad|\psi|^{(2+\a,1+\a/2)}_{e_{0,1/2}\times[0,T_{max}]}\leq K |m_0/\phi(\cdot,0)|^{(2+\a)}_{\G}\\[12pt]
 (iii)&\textrm{If $m_0$   fulfills \eqref{cmpv0} and \eqref{FPallaVB}:}\qquad|\psi|^{(2+\a,1+\a/2)}_{\G\times[0,T_{max}]}\leq K |m_0/\phi(\cdot,0)|^{(2+\a)}_{\G}.
\end{array}
\end{equation*}
where $K$ is a constant independent of $m_0$ and $c_T$.
\end{Theorem}
\begin{Proof}
Proposition \ref{exisheat} ensures all the part of the statement concerning the function $\phi$.
Invoking Proposition \ref{exisFP0} (respectively, Proposition \ref{exisFP}-($i$) and -($ii$)), by the regularity and the lower bound of $\phi$, we deduce the part of the statement concerning the function $\psi$ in point ($i$) (respectively, in point ($ii$) and in point ($iii$)).
\end{Proof}
We   also have existence and uniqueness for the solution  to \eqref{MFGtoy}:
\begin{Corollary}
Under the hypotheses of Theorem \ref{thm:MFG},  there exist a unique classical solution to the MFG system \eqref{MFGtoy}.
\end{Corollary}
Being a straightforward consequence of the previous theorem, the proof of this result is omitted.

%%%%%%%%%%%%
%          %
%%%%%%%%%%%%

\section{The Mean Field Game result}\label{S4}
We   prove the existence of a   starting time $T$ consistent with the corresponding flux of participants $\pd_x m$. To this end we show that the map from $[t_0,T_{max}]$ into itself, defined by the  scheme \eqref{scheme} is continuous and therefore it admits a fixed point by the Brouwer's Theorem. For simplicity, we shall recast it in terms of couple $(\phi,\psi)$ solution of \eqref{MFGchange}.
Consider the function $\Psi: [t_0,T_{max}]\to [t_0,T_{max}]$ defined as
\begin{equation}\label{defPsi}
T\rightarrow c_T\rightarrow \phi\rightarrow \psi\rightarrow T^*=:\Psi(T)
\end{equation}
where $T^*$ is defined as in \eqref{quorum} with
\begin{equation}\label{Fheat}
F(s)=\int_0^s e^{c_T(r)}\, \pd_x \psi(v_0,r)\, dr=:\int_0^s\tilde\psi_T(r)\, dr.
\end{equation}
%here, the initial data of system \eqref{MFGchange} have been used (and the regularity of $\phi$ and $\psi$ stated in Theorem \ref{thm:MFG} as well).\\
In this section we assume the hypotheses of Theorem~\eqref{thm:MFG} and that
 the map
\begin{equation}\label{Creg}
%\begin{split}
%c_T:&[0,T_{max}]\rightarrow C^{(1+\a/2)}([0,T_{max}])\\
 T\in [0,T_{max}]\mapsto c_T \in C^{(1+\a/2)}([0,T_{max}])
%\end{split}
\end{equation}
is continuous.
A crucial  step to prove the existence of the mean field $T$  is to establish some bounds for $\pd_x\psi(v_0,\cdot)$. In order to get such an estimate, we consider in the next Lemma two complementary assumptions.

\begin{Lemma}\label{lemma:GLL2.5}
Let $(\phi,\psi)$ be the solution to system \eqref{MFGchange}.
\begin{itemize}
\item[(a)] If
\begin{equation}\label{cmpv1a}
\pd_x m_0(v_0)>0,
\end{equation}
then, there exists a value $\e>0$, independent of $T$, such that
\begin{equation*}
|\pd_x \psi(v_0,t)|>\e \qquad \forall t\in [0,T_{max}].
\end{equation*}
\item[(b)] If $m_0$ fulfills \eqref{cmpv0},
%If
%\begin{equation}\label{cmpv1b}
%\pd_x  m_0(v_0)=\pd^2_x m_0(v_0)=0.
%\end{equation}
then there holds:
\begin{equation*}
\pd_x\psi(v_0,t)>0\qquad \forall t\in(0,T_{max}].
\end{equation*}
In particular,
there exists a constant $\e_T$ such that
\begin{equation*}
|\pd_x \psi(v_0,t)|>\e_T  \qquad \forall t\in [t_0,T_{max}].
\end{equation*}
\end{itemize}
\end{Lemma}
\begin{Proof}
\texttt{($a$).}
%Since  the point $v_0$ is a boundary vertex  there exists a unique edge, say $e_0$, which is incident to it. Without any loss of generality, we denote $v_1$ the other endpoint of $e_0$ and we assume that the parametrization of $e_0$ fulfills: $\pi_0(0)=v_0$ and $\pi(l_0)=v_1$.
Owing to \eqref{ipminime}, the function $m_0$ satisfies: $m_0(v_0)=0$ and $\pd_x m_0(v_0)>0$. Moreover, Proposition \ref{exisheat} ensures that $|\frac{m_0}{\phi(\cdot,0)}|^{(2+\a)}_{[0,l_0]}$ is bounded independently of $T$.
We infer that there exist $\xi_0\in(0,l_0)$ and a sufficiently small $a>0$ such that, for every $T\in[0,T_{max}]$ there holds
\begin{equation*}
\frac{m_0(x)}{\phi(x,0)}\geq a \sin\left(\frac{x\pi}{\xi_0}\right)\qquad \forall x\in[0,\xi_0].
\end{equation*}
One can easily check that the function
\[v(x,t):=a e^{bt} \sin(x\pi/\xi_0),\qquad\textrm{with}\qquad b:=-\pi^2/\xi_0^2\]
solves the initial-boundary value problem
\[\left\{
\begin{array}{ll}
\pd_t v- \pd_x^2 v=0&\qquad (x,t)\in (0,\xi_0)\times (0,T_{max})\\
v(0,t)=v(\xi_0,t)=0&\qquad t\in (0,T_{max})\\
v(x,0)=a \sin(x\pi/\xi_0)&\qquad x\in (0,\xi_0)
\end{array}\right.\]
while the function $\psi$ is a supersolution to this problem.
By the standard comparison principle, we infer: $\psi\geq v$ in $[0,\xi_0]\times [0,T_{max}]$.
Since $\psi(0,\cdot)=v(0,\cdot)$ on $[0,T_{max}]$, we get $\pd_x\psi(0,t)\geq \pd_x v(0,t)=a e^{bt} \pi/\xi_0$.
In particular, we deduce
\[
|\pd_x\psi(0,t)|\geq a e^{bT_{max}} \pi/\xi_0\qquad \forall t\in[0,T_{max}]
\]
where all the constants are independent of $T$.
%
%fine dim punto (a)
%

\texttt{($b$).}
Being nonnegative, the function $\psi$ attains a global minimum at each point $(v_0,t)$ with $t\in(0,T_{max}]$. The Hopf Lemma prevents that  $\pd_x\psi(v_0,t_0)\leq 0$ in these points. Hence, there holds: $\pd_x\psi(v_0,t)>0$ in $(0,T_{max}]$. The second part of the statement follows by continuity.
\end{Proof}

%\begin{Remark}
%The proof of the previous lemma also provides an almost explicit formula for the value $\e$. Without obtaining such an estimate, the statement can be established arguing as in \cite[Lemma 2.5]{gll}: for each $T$, the function $\psi$ attains its global minimum at every point $(v_0,t)$, $t\in[0,T_{max}]$. The Hopf Lemma prevents that $\pd_x\psi (v_0,t)=0$ for any $t\in[0,T_{max}]$. Since $\psi$ is regular and $\pd_x\psi (v_0,0)\ne 0$, by continuity we infer that $|\pd_x \psi(v_0,\cdot)|>0$ in  $[0,T_{max}]$. Finally, since the solution $(\phi,\psi)$ to system \eqref{MFGchange} depends continuously on the features of the problem, we get the statement of Lemma \ref{lemma:GLL2.5}.
%\end{Remark}

We shall establish the existence of a fixed point provided that $m_0$ fulfills either \eqref{cmpv1a} or  \eqref{cmpv0}. We cope with these two cases separately in the next two statements.
\begin{Theorem}\label{MFGheat1}
Assume the hypotheses of Theorem~\eqref{thm:MFG}-($i$) and inequality~\eqref{cmpv1a}. Then the map $\Psi:[0,T_{max}]\to [0,T_{max}]$ defined by \eqref{defPsi} admits a fixed point.
\end{Theorem}
\begin{Proof}
%To prove the result it is expedient to introduce the function \[\Psi^*:\, C^0([0,T_{max}])\times C^1(\G\times[0,T_{max}])\rightarrow [0,T_{max}]\] defined by $\Psi^*(c_T,\psi):=T^*$ where $T^*$ is the value defined in \eqref{quorum} with $F$ given by \eqref{Fheat}.
%The function $\Psi^*$ is well defined because of the uniqueness of the solution to  \eqref{MFGchange}.\\
%Denote by $(\phi_T,\psi_T)$ the solution   to system \eqref{MFGchange} corresponding cost $c_T$.
%We prove tat the function $\Psi^*$ is continuous if it is restricted to the set
%\[{\cal A}:=\{(c_T,\psi_T) \mid T\in[0,T_{max}]\}\subset C^{(1+\a/2)}([0,T_{max}])\times C^{(2+\a,1+\a/2)}(\G\times [0,T_{max}])\]
%(the inclusion is a consequence of Theorem \ref{thm:MFGc}).
We shall follow the arguments of \cite[Lemma 2.6]{gll}. In order to apply the  Brouwer fixed point Theorem, we need to prove that the function $\Psi$ defined in \eqref{defPsi} is continuous.
We consider two admissible flows $\tilde \psi_{T_1},\tilde \psi_{T_2}$ (see equation \eqref{Fheat} for their definition)
% namely (for $i=1,2$) $\tilde\psi_i:= e^{\frac{c_{T_i}(r)}{\s^2}}\, \pd_x \psi_i(v_0,r)$ where $(\phi_i,\psi_i)$ solves system \eqref{MFGchange} with $c_T=c_{T_i}$.We set $T^*_i:=\Psi^*(c_{T_i},\psi_i)$ for $i=1,2$
and, without any loss of generality, we assume $\Psi(T_1)\leq \Psi(T_2)$. If $\Psi(T_1), \Psi(T_2)\in(t_0,T_{max})$, we have
\begin{equation*}
0=\int_0^{\Psi(T_1)}\tilde \psi_{T_1}(t)dt-\int_0^{\Psi(T_2)}\tilde \psi_{T_2}(t)dt=
\int_0^{\Psi(T_1)}(\tilde \psi_{T_1}(t)-\tilde\psi_{T_2}(t))dt-\int_{\Psi(T_1)}^{\Psi(T_2)}\tilde \psi_{T_2}(t)dt
\end{equation*}
(where the first equality is due to the fact that both integrals are equal to $\th$). Taking into account Lemma \ref{lemma:GLL2.5}-($a$), we obtain
\begin{equation*}
\e (\Psi(T_2)-\Psi(T_1))\leq \int_0^{\Psi(T_1)}(\tilde \psi_{T_1}(t)-\tilde\psi_{T_2}(t))dt\leq |\tilde \psi_{T_1}-\tilde\psi_{T_2}|_{L^1(0,T_{max})}.
\end{equation*}
The estimates in Theorem \ref{thm:MFG}-($i$) and the trace theorem (for instance, see \cite[Theorem II.2.3]{lsu}) yield
\begin{equation*}
\Psi(T_2)-\Psi(T_1)\leq \textrm{const.}\, |c_{T_1}-c_{T_2}|^{(1+\a/2)}_{[0,T_{max}]}.
\end{equation*}
Taking into account assumption \eqref{Creg}, we obtain that in this case the function $\Psi$ is continuous.

When $\Psi(T_1)=t_0$ (respectively, $\Psi(T_2)=T_{max}$), we have
\[
\int_0^{\Psi(T_1)}\tilde \psi_{T_1}-\int_0^{\Psi(T_2)}\tilde \psi_{T_2}\geq 0;
\]
indeed, either  $\tilde \psi_{T_1}$ is a flux which reaches $\th$ at most at time $\Psi(T_1)$ or $\tilde \psi_{T_2}$ is a flux which does not reach the value $\th$ before time $T_{max}$ ; in other words, the former integral is $\geq \th$ (respectively, the latter one is $\leq \th$).
Hence we can conclude by the same arguments as before.
Therefore, the continuity of $\Psi$ is achieved.
\end{Proof}
%\begin{Remark}
%Since $\pd_x\psi(v_0,t_0)>0$, the second part of the statement can also be obtained arguing as in \cite[Lemma 2.5]{gll} or as follows: observe that, for suitable $a$ and $\xi_0$ and $b:=-\nu\pi^2/\xi_0^2$, the function $a e^{b(t-t_0)}\sin(x\pi/\xi_0)$ is a subsolution to the Fokker-Planck equation in $e_{0,\xi_0/l_0}\times(t_0,T_{max})$.
%\end{Remark}

\begin{Theorem}\label{MFGheat}
Assume the hypotheses of Theorem \ref{thm:MFG}-($ii$). Then the map $\Psi$ defined by \eqref{defPsi} admits a fixed point.
\end{Theorem}
\begin{Proof}
We shall argue adapting the arguments of Theorem \ref{MFGheat1}: hence, our purpose is to prove that $\Psi$ is continuous on $[t_0,T_{max}]$.
To this end, let us fix $T\in[t_0,T_{max}]$.
%denote $c_T$, $(\phi_T,\psi_T)$ and $\e_T$ respectively the corresponding cost, the solution to system \eqref{MFGchange} and the value $\e_{T}>0$ defined in Lemma \eqref{lemma:GLL2.5}-($b$). We want to prove that $\Psi$ is continuous in $T$.
For every $T_1\in[t_0,T_{max}]$ such that $\Psi(T)=\Psi(T_1)$, there is nothing to prove.
%We set
%\begin{equation}\label{def:tildepsi}
%\tilde\psi_T(s):=\exp\{c_T(s)/\s^2\}\pd_x\psi_T(v_0,s)\qquad \forall s\in[0,T_{max}]
%\end{equation}
%(note that the right hand side coincides with the integrand in formula \eqref{Fheat}).
We split the arguments according to the fact that $\Psi(T)$ belongs to $(t_0,T_{max})$, to $\{t_0\}$ or to $\{T_{max}\}$.

\texttt{Case1: $\Psi(T)\in(t_0,T_{max})$.}
% Observe that, by Theorem \eqref{thm:MFG}, there exists a neighborhood ${\cal V}$ of $T$ such that $\Psi({\cal V})\subset(t_0,T_{max})$.
Consider $T_1\in[t_0,T_{max}]$ with $\psi(T_1)<\Psi(T)$; set
\begin{equation}\label{def:tau}
\t:=\inf\{t\in(0,T_{max})\mid \int_0^t\tilde\psi_{T_1}=\th\}
\end{equation}
and observe that $\Psi(T_1)=\max\{t_0,\tau\}$.
Then, we have
\begin{equation*}
0=\int_0^{\t}\tilde\psi_{T_1}-\int_0^{\Psi^(T)}\tilde\psi_{T}=
\int_0^{\t}\left(\tilde\psi_{T_1}-\tilde\psi_{T}\right)-
\int_{\t}^{\Psi^(T)}\tilde\psi_{T}
\end{equation*}
(the first equality is due to the fact that both the integrals are equal to $\th$). By Lemma \ref{lemma:GLL2.5}-($b$), we infer
\begin{equation*}
\e_T\left(\Psi(T)-\Psi(T_1)\right)\leq \int_{\Psi^(T_1)}^{\Psi^(T)}\tilde\psi_{T}
\leq \int_{\t}^{\Psi^(T)}\tilde\psi_{T}
=\int_0^{\t}\left(\tilde\psi_{T_1}-\tilde\psi_{T}\right)\leq |\tilde \psi_{T_1}-\tilde\psi_{T}|_{L^1(0,T_{max})}
\end{equation*}
Arguing as before, we deduce that there exists a constant $\tilde K$ (depending on $T$) such that
\begin{equation}\label{cont1}
\Psi(T)-\Psi(T_1)\leq \tilde K |T_1-T|.
\end{equation}
%Then, we have
%\begin{equation*}
%0=\int_0^{\Psi^(T_1)}\tilde\psi_{T_1}-\int_0^{\Psi^(T)}\tilde\psi_{T}=\int_0^{\Psi^(T_1)}\left(\tilde\psi_{T_1}-\tilde\psi_{T}\right)-\int_{\Psi^(T_1)}^{\Psi^(T)}\tilde\psi_{T}
%\end{equation*}
%(the first equality is due to the fact that both the integrals are equal to $\th$). By Lemma \ref{lemma:GLL2.5bis}, we infer
%\begin{equation*}
%\e_T\left(\Psi(T)-\Psi(T_1)\right)\leq \int_{\Psi^(T_1)}^{\Psi^(T)}\tilde\psi_{T}=\int_0^{\Psi^(T_1)}\left(\tilde\psi_{T_1}-\tilde\psi_{T}\right)\leq  T_{max}\|\tilde\psi_{T_1}-\tilde\psi_T\|_\infty.
%\end{equation*}
%By Theorem \ref{thm:MFG}, there exists a constant $\tilde K$ (depending on $T$) such that
%\begin{equation}\label{cont1}
%|\Psi(T)-\Psi(T_1)|\leq \tilde K |T_1-T|.
%\end{equation}

Consider now a point $T_1\in[t_0,T_{max}]$ with $\psi(T_1)>\Psi(T)$. Then, we have
\begin{equation*}
0\leq\int_0^{\Psi^(T)}\tilde\psi_{T}-\int_0^{\Psi^(T_1)}\tilde\psi_{T_1}=\int_0^{\Psi^(T_1)}\left(\tilde\psi_{T}-\tilde\psi_{T_1}\right)+
\int_{\Psi^(T_1)}^{\Psi^(T)}\tilde\psi_{T}
\end{equation*}
where the inequality is due to the fact that the first integral is equal to $\th$ while the second one is less or equal to $\th$.
Again by Lemma \ref{lemma:GLL2.5}-($b$), we infer
\begin{equation*}
\e_T\left(\Psi(T_1)-\Psi(T)\right)\leq \int_{\Psi^(T)}^{\Psi^(T_1)}\tilde\psi_{T}\leq\int_0^{\Psi^(T_1)}\left(\tilde\psi_{T}-\tilde\psi_{T_1}\right)\leq
|\tilde \psi_{T_1}-\tilde\psi_{T}|_{L^1(0,T_{max})}.
\end{equation*}
Arguing as before, for some constant $\tilde K'$ (depending on $T$), we get
\begin{equation*}
\Psi(T_1)-\Psi(T)\leq \tilde K' |T_1-T|.
\end{equation*}
By this relation and \eqref{cont1}, the proof of the continuity of $\Psi$ in $T$ is accomplished.

\texttt{Case2: $\Psi(T)=T_{max}$.} For $T_1\in[t_0,T_{max}]$ with $\Psi(T_1)=T_{max}$, there is nothing to prove; hence, without any loss of generality, we assume that $\Psi (T_1)<T_{max}$.
We have
\begin{equation*}
0\leq\int_0^{\Psi^(T_1)}\tilde\psi_{T_1}-\int_0^{T_{max}}\tilde\psi_{T}=\int_0^{\Psi^(T_1)}\left(\tilde\psi_{T_1}-\tilde\psi_{T}\right)-
\int_{\Psi^(T_1)}^{T_{max}}\tilde\psi_{T}
\end{equation*}
Arguing as before, we accomplish the proof in this case.

\texttt{Case3: $\Psi(T)=t_0$.} For $T_1\in[t_0,T_{max}]$ with $\Psi(T_1)=t_0$, there is nothing to prove; hence, without any loss of generality, we assume that $\Psi (T_1)>t_0$.
We have
\begin{equation*}
0\leq \int_0^{t_0}\tilde\psi_{T}-\int_0^{\Psi^(T_1)}\tilde\psi_{T_1}=\int_0^{\Psi^(T_1)}\left(\tilde\psi_{T}-\tilde\psi_{T_1}\right)+
\int_{\Psi^(T_1)}^{t_0}\tilde\psi_{T}.
\end{equation*}
%where the inequality is due to the following fact: $\tilde \psi_T$ is a flux such that the value $\th$ is reached at most at time $t_0$ (i.e.: $\int_0^{t_0}\tilde\psi_T\geq \th$). Lemma \ref{lemma:GLL2.5bis} ensures
%\begin{equation*}
%\e_T\left(\Psi(T_1)-t_0\right)\leq \int_{t_0}^{\Psi^(T)}\tilde\psi_{T}=\int_0^{\Psi^(T_1)}\left(\tilde\psi_{T}-\tilde\psi_{T_1}\right)\leq T_{max}\|\tilde\psi_{T_1}-\tilde\psi_T\|_\infty.
%\end{equation*}
By the same arguments as those used before, we accomplish the proof.
\end{Proof}

%Let us state the consequences of the previous results on  Mean Field problem discussed in Section  \ref{toymodel}.
\begin{Corollary}
Under the hypotheses of either Theorem \ref{MFGheat1} or Theorem \ref{MFGheat}, there exists a value $T$ which is coherent with the expectation of the participants to the meeting.
\end{Corollary}

We conclude with a uniqueness result for the fixed point under some monotonicity condition on the cost $c_T$.
\begin{Proposition}
Assume that the cost $c_T$ does not depend on the term $c_2$, then the map $\Psi$ defined by \eqref{defPsi} admits a unique fixed point.
\end{Proposition}
\begin{Proof}
Existence of a fixed point is proved in either Theorem~\ref{MFGheat1} or Theorem~\ref{MFGheat}. Assume by contradiction that there exist $T_1, T_2\in [0, T_{max}]$ with $T_1>T_2$ such that $T_i=\Psi(T_i)$. Let $c_{T_i}$ and  $(\phi_i, \psi_i)$ be the costs and the  solutions of \eqref{MFGchange}  corresponding to $T_i$, $i=1,2$.
Then, $(\phi,\psi):=(\phi_1-\phi_2,\psi_1 -\psi_2)$ satisfies \eqref{MFGchange} with $m_0/\phi(\cdot,0)$, $e^{c_T(T_{max})}$ and $e^{c_T(\cdot)}$ replaced respectively by $m_0/\phi_1(\cdot,0)-m_0/\phi_2(x,0)$, $e^{c_{T_1}(T_{max})}-e^{c_{T_2}(T_{max})}$ and $e^{c_{T-1}(\cdot)}-e^{c_{T_2}(\cdot)}$.
%\begin{equation}\label{uniq1}
%  \left\{
%    \begin{array}{ll}
%     -\pd_t \phi-\nu\pd^2_x\phi=0\quad &(x,s)\in \G\times (0,T_{max}),\\[6pt]
%     \pd_t \psi-\nu\pd^2_x\psi=0       &(x,s)\in \G\times (0,T_{max})\\[6pt]
%     \psi(x,0)=\frac{m_0(x)}{\phi_1(x,0)}-\frac{m_0(x)}{\phi_2(x,0)},\, \phi(x,T_{max})=e^{ \frac{c_{T_1}(T_{max})}{\s^2}}-e^{ \frac{c_{T_2}(T_{max})}{\s^2}}   &x \in\G\\[6pt]
%     \psi(v_0,s)=0,\, \phi(v_0,s)=e^{ \frac{c_{T_1}(s)}{\s^2}}-e^{ \frac{c_{T_2}(s)}{\s^2}}& s\in[0, T_{max}]
%   \end{array}
% \right.
%\end{equation}
We have
\begin{align*}
    0=&\int_0^{T_{max}}\int_\G [-\pd_t \phi-\pd^2_x\phi]\psi dx\,dt=
       \int_0^{T_{max}}\int_\G   [\pd_t \psi\,  \phi+\pd_x\phi\,\pd_x\psi]  dx\,dt \\
       &-\int_\G\big[\psi(x,\cdot)\phi(x,\cdot)\big]_0^{T_{max}}dx-\sum_{i\in I}\sum_{j\in Inc_i}\int_0^{T_{max}}(-a_{ij}) \pd_j \phi(v_i,t) \psi(v_i,t) dt
\end{align*}
(the term  $-a_{ij}$ takes into account the orientation of the edge $e_j$).
Similarly
\begin{align*}
    0=&\int_0^{T_{max}}\int_\G [\pd_t \psi-\pd^2_x\psi]\phi dx\,dt=\int_0^{T_{max}}\int_\G   [\pd_t \psi\,\phi+\pd_x\phi\pd_\psi ] dx\,dt \\
    &-\sum_{i\in I}\sum_{j\in Inc_i}\int_0^{T_{max}}(-a_{ij})\pd_j \psi(v_i,t) \phi(v_i,t) dt
\end{align*}
Subtracting the previous inequality and  using the transition conditions at the internal nodes   we get
\begin{align*}
0=&  \int_\G\left(\frac{m_0(x)}{\phi_1(x,0)}-\frac{m_0(x)}{\phi_2(x,0)}\right)(\phi_1(x,0)-\phi_2(x,0))dx\\
&-( e^{c_{T_1}(T_{max})}-e^{c_{T_2}(T_{max})} )\int_\G(\psi_1(x,T_{max})-\psi_2(x,T_{max}))dx\\
 &\int_0^{T_{max}}(e^{c_{T_1}(t)}-e^{c_{T_2}(t)})(\pd_{0}\psi_1(v_0,t)-\pd_{0}\psi_2(v_0,t)) dt
\end{align*}
(recall that $e_{0}$ is the unique arc incident to $v_0$ parameterized in such a way that $v_0$ is the initial point).\\
The first term in the previous inequality is negative.
By the assumption on $c_T$, the map    $T\mapsto c_T$ is increasing in $T$  and $c_{T_1}(T_{max})=c_{T_2}(T_{max})$. Hence the second term is null. Moreover, since  $T_1>T_2$ and therefore $c_{T_1} > c_{T_2}$ on $[0, T_{max}]$, we have $\phi_1\ge \phi_2$, hence  $\psi_1\le \psi_2$ and,  by $\psi_i(v_0,t)=0$ for $i=1,2$,  $\pd_{0}\psi_1(v_0,t) \le \pd_{0} \psi_2(v_0,t)$. It follows that also the third term is negative, hence $\phi_1(x,0)=\phi_2(x,0)$ for $x\in\G$ and therefore a contradiction to $c_{T_1}> c_{T_2}$.
\end{Proof}
%%%%%%%%%%%%%%%%%%
%                %
%%%%%%%%%%%%%%%%%%
\section{Numerical simulation}\label{S5}
In this section we propose a numerical method to compute the mean field $T$. The scheme is based on a finite difference approximation of the system \eqref{MFGchange} with an iterative procedure to solve the fixed point map  \eqref{defPsi}.\\
On each interval  $[0,l_j]$, $j\in J$,  it is defined an uniform  partition $y_k=k h_j$ with space step $h_j=\frac{l_j}{M_j}$ and $k=0,\dots,M_j$. In this way   a spatial grid $\mathcal{G}(\Gamma)=\{x_{j,k}=\pi_j(y_k),\; j\in J, \; k=0,\dots,M_j\}$ is  defined on the network $\Gamma$. A time step  $\D t$ is also introduced  to obtain  a uniform grid $t_n=n\Delta t$, $n=0,1,\dots,N_{max}$ with $N_{max}=[ T_{max}/\Delta t ]$, on the time interval $[0,T_{max}]$.\\%(where $\lceil\cdot \rceil$ denotes the upper integer part) .\\
%We  denote   the value  of a function $f:\mathcal{G}(\Gamma)\to\R$,  on the grid  node $x_{j,k}$,  by  $f_{j,k}$.   \\
We will approximate the solution $(\phi,\psi)$ of \eqref{MFGchange} by  two sequences  $\{\phi^n\}_n$ and $\{\psi^n\}_n$, where,  for each $n=0,\dots,N_{max}$, $\phi^n,\psi^n:\mathcal{G}(\Gamma)\to\R$ and   $\phi^n_{j,k}\simeq \phi(x_{j,k},t_n)$, $\psi^n_{j,k}\simeq \psi(x_{j,k},t_n)$.
The discrete functions  $\{\phi^n\}_n$ and $\{\psi^n\}_n$ are computed by the  following  forward-backward explicit finite difference scheme:
 %for $k=1,\dots,M-1$ and $ j\in J$,
\begin{equation}{\label{FDscheme}}
\begin{cases}
\phi^{n}_{j,k}=\phi^{n+1}_{j,k}+\displaystyle{\frac{\D t}{h_j^2}}\left(\phi^{n+1}_{j,{k+1}}-2\phi^{n+1}_{j,k}+\phi^{n+1}_{j,k-1}\right),
&\; n=N_{max}-1,\dots,0\\
\psi^{n+1}_{j,k}=\psi^{n}_{j,k}+\displaystyle{\frac{\D t}{ h_j^2}}\left(\psi^{n}_{j,k+1}-2\psi^{n}_{j,k}+\psi^{n}_{j,k-1}\right),
&\; n=0,\dots,N_{max}-1\\
{\mbox{ for }} k=1,\dots,M_j-1\,{\mbox {and}} \; j\in J.
\end{cases}
\end{equation}
At each time iteration $n$, to compute $\{\phi^n\}_n$ and $\{\psi^n\}_n$ it is necessary to fix the values of these functions at the boundary of the arcs $e_j$, $j\in J$, i.e. at the transition vertices $v_i$, $i\in I_T$. We define an approximation of the Kirchhoff's condition  which together with the continuity condition across the vertices will give the $\#(Inc_i)$ conditions necessary to determine in a unique way the value of the functions $\phi^n$ and $\psi^n$ at $v_i$.\\
%These sequence verify a scheme that we are going to derive.  \\
%and  approximate the solution $(\phi,\psi)$ of \eqref{MFGchangefull} on the nodes, $\phi^n_{j,k}\simeq \phi(x^j_k,t_n)$ and  $\psi^n_{j,k}\simeq \psi(x^j_k,t_n)$.\\
%Before introducing the discrete system the functions $\{\phi^n, \psi^n\}_n$  have to  solve,
We introduce two sets of indices $Inc^+_i=\{j\in J\;| a_{ij}=1 \}$ and $Inc^-_i=\{j\in J\;| a_{ij}=-1 \}$. Moreover we denote by $\phi^n(v_i)$, $\psi^n(v_i)$  the values of the functions $\phi^n$, $\psi^n$  at  $v_i\in V$. If $j\in Inc^+_{i}$, then $\phi(\pi_j(y_0),t_n)\simeq \phi^n_{j,0}=\phi^n(v_i)$ while  if $j\in Inc^-_{i}$, then $\phi(\pi_j(y_{M_j}),t_n)\simeq \phi^n_{j,M_j}=\phi^n(v_i)$.
We   define the following finite differences approximations of the  derivatives at $v_i$ along an edge $e_j$:
% \eqref{Darc1}-\eqref{Darc2} as:
\begin{align*}
 \partial_{j}\phi(v_i,t_n)\simeq \frac{1}{h_j}\left(\phi^n_{j,1}-\phi^n(v_i)\right),\;
&\partial_{j}\psi(v_i,t_n)\simeq \frac{1}{h_j}\left(\psi^n_{j,1}-\psi^n(v_i)\right) \; &j\in Inc^+_{i},\\
 \partial_{j}\phi(v_i,t_n)\simeq \frac{1}{h_j}\left(\phi^n_{j,M_j-1}-\phi^n(v_i)\right),\;
& \partial_{j}\psi(v_i,t_n)\simeq \frac{1}{h_j}\left(\psi^n_{j,M_j-1}- \psi^n(v_i)\right) \;&j\in Inc^-_{i}.
\end{align*}
We  rewrite the transition conditions in \eqref{MFGchange} as
\begin{eqnarray}
\label{TC2a}\sum_{j \in Inc^+_{i}} \partial_{j}\phi(v_i,s)-\sum_{j \in Inc^-_{i}}\partial_{j}\phi(v_i,s)=0,\quad \\
\label{TC2b}\sum_{j \in Inc^+_{i}}\partial_{j}\psi(v_i,s)-\sum_{j \in Inc^-_{i}}\partial_{j}\psi(v_i,s)=0\quad \end{eqnarray}
%with $(v_i,s) \in \Gamma\times (0,T_{max}).$
%We consider  the following  approximation of  the derivative along an edge $e_j$ :
%\begin{eqnarray*}
%&&\partial_{x,j}\phi(\pi_j(y_0),s)\simeq \frac{1}{h_j}(\phi(\pi_j(y_{1}),s)-\phi(\pi_j(y_0),s))\label{Darc1}\\
%&&\partial_{x,j}\phi(\pi_j(y_M),s)\simeq\frac{1}{h_j}( \phi(\pi_j(y_{M-1}),s)-\phi(\pi_j(y_M),s)).\label{Darc2}\\
%&&\partial_{x,j}\phi(\pi_j(y_k),s)\simeq \frac{1}{2h_j}(\phi(\pi_j(y_{k+1}),s)-\phi(\pi_j(y_{k-1}),s))\; k=1,M-1.\nonumber
%\end{eqnarray*}
and we consider the following finite difference approximation
\begin{eqnarray}
&\sum_{j\in Inc^+_{i}} \frac{1}{h_j}(\phi^n_{j,1}-\phi^n(v_i))-\sum_{j\in Inc^-_{i}}\frac{1}{h_j}(\phi^n(v_i)-\phi^n_{j,M_j-1})=0,\quad \quad \label{Kdis1}\label{contcond2a}\\
&\sum_{j\in Inc^+_{i}} \frac{1}{h_j}(\psi^n_{j,1}-\psi^n(v_i))-\sum_{j\in Inc^-_{i}}\frac{1}{h_j}(\psi^n(v_i)-\psi^n_{j,M_j-1})=0.\quad\quad  \label{Kdis2}
\label{contcond2b}
\end{eqnarray}
%Let us explicit  with respect $\phi^n(v_i) $ and $\psi^n(v_i)$
%\begin{eqnarray}
%\phi^n(v_i) =\frac{\sum_{j\in Inc^+_{i}}\frac{a_{i,j}\rho_{i,j}}{h_j}\phi^n_{1,j}+\sum_{j\in Inc^-_{i}}\frac{a_{i,j}\rho_{i,j}}{h_j}\phi^n_{M-1,j}}{\sum_{j\in Inc_{i}}\frac{a_{i,j}\rho_{i,j}}{h_j}},\label{contcond2a}\\
%\psi^n(v_i) =\frac{\sum_{j\in Inc^+_{i}}\frac{a_{i,j}\rho_{i,j}}{h_j}\psi^n_{1,j}+\sum_{j\in Inc^-_{i}}\frac{a_{i,j}\rho_{i,j}}{h_j}\psi^n_{M-1,j}}{\sum_{j\in Inc_{i}}\frac{a_{i,j}\rho_{i,j}}{h_j}}.\label{contcond2b}
%\end{eqnarray}
Given a discrete function $f:\mathcal{G}(\Gamma)\to\R$, we consider a continuous piecewise linear reconstruction $I[f]:\Gamma\to \mathbb{R}$  such that $I[f]\mid_{(x_{j,k},x_{j,k+1})}$ is linear for all $j\in J$ and $k=0,\dots,M_j-1$ and $I[f](x_{j,k})=f_{j,k}$. To guarantee the continuity on $\Gamma$ of  the linear interpolation $I[\cdot]$ applied to the discrete function  $\phi^n$ and $\psi^n$, we need to impose the  following  continuity conditions:
\begin{align}
% \mbox {for all  $i\in I_T$}
&\phi^n_{j,0}=\phi^n(v_i), \quad \psi^n_{j,0}=\psi^n(v_i)  \quad  &\text{if $i\in I_T$, $j \in Inc^+_i$},  \label{contconda}\\
&\phi^n_{j,M_j}=\phi^n(v_i), \quad \psi^n_{j,M_j}=\psi^n(v_i) \quad & \text{if $i\in I_T$, $j \in Inc^-_i$}.  \label{contcondb}
\end{align}
At each time step $t_n$, the  $\# (Inc_i)-1$ conditions given by \eqref{contconda}-\eqref{contcondb}   coupled with \eqref{contcond2a}-\eqref{contcond2b} give $\#(Inc_i)$ relations which uniquely determine    $\phi^n(v_i)$ and     $\psi^n(v_i)$.\\
Summarizing, we approximate  \eqref{MFGchange} by computing the couple of discrete functions $\{(\phi^n,\psi^n)\}_n$
 which solve the  finite difference scheme \eqref{FDscheme}
together with
\begin{description}
\item[$i)$] the   conditions \eqref{contcond2a}-%\eqref{contcond2b} and respectively \eqref{contconda}-
\eqref{contcondb}  at the vertices $v_i\in\G_T$;
\item[$ii)$]  the boundary condition
\begin{equation*}\label{BC}
\phi^n(v_0)=e^{c_T(t_n)}\quad \psi^n(v_0)=0 \quad n=0,\dots ,N_{max};
\end{equation*}
\item[$iii)$] the initial and terminal conditions:
\begin{equation*}
 \phi^{N_{max}}_{j,k}=e^{c_T(T_{max})}, \quad\psi^0_{j,k}=\frac{m_0(x_{j,k})}{\phi^0_{j,k}}, \quad k=0,\dots,M_j-1,\,j\in J.
\end{equation*}
\end{description}
Defined a function $\{\psi^n\}_n$ by means of the previous scheme, we consider
the following  approximation of the cumulative distribution
 \eqref{Fheat}
\begin{equation}{\label{CostTime}}
\tilde F(t_n)=\frac{\D t}{h_0}\sum_{k=0}^n e^{c_T(k\D t)}\psi^k_{0,1},
\end{equation}
where      $e_0$ denotes the edge incident   $v_0$ with $\pi_0(0)=v_0$ and by the boundary condition $\psi^k_{0,0}=\psi^k(v_0)=0$.\\
To approximate the fixed point of the map  $\Psi$ defined in  \eqref{defPsi} we apply the following Algorithm \ref{algo}. Given an initial guess $T_1$ and denoted by $T_2$ an initial value   to enter the loop  and by    $\tau$ as threshold for the stopping criteria, we consider \\

%\IncMargin{1em}
\begin{algorithm}[H]\label{algo}
%\SetAlgoLined
\KwData{  initial guess  $T_1$, $T_2$,  threshold value $\tau$}
\KwResult{approximated mean field $T_2$}
\While{$|T_1-T_2|> \tau$}{
 set $T_1\gets T_2$\;
 %\label{step2}
 solve \eqref{FDscheme} with $T=T_1$ and conditions  $i),ii),iii)$\;
 compute
\begin{equation*}{\label{minCostTime}}
T_{N^*}=\min\{n\D t, n=0,\dots,N_{max}| \tilde F(t_n)>\theta\}
\end{equation*}
%\State set $T_2=N_2 \D t$;
\eIf{$T_{N^*}<t_0$}{
    set $T_2\gets t_0$;
      }{
      set $T_2\gets T_{N^*}$; }
}
\caption{Fixed Point Iterations}
\end{algorithm}
%\DecMargin{1em}
%where $\tau$ is a given threshold for the stopping criteria the fixed point algorithm

\subsection{Example 1: a simple graph}
\begin{figure}[ht!!]
\begin{center}
\epsfig{figure=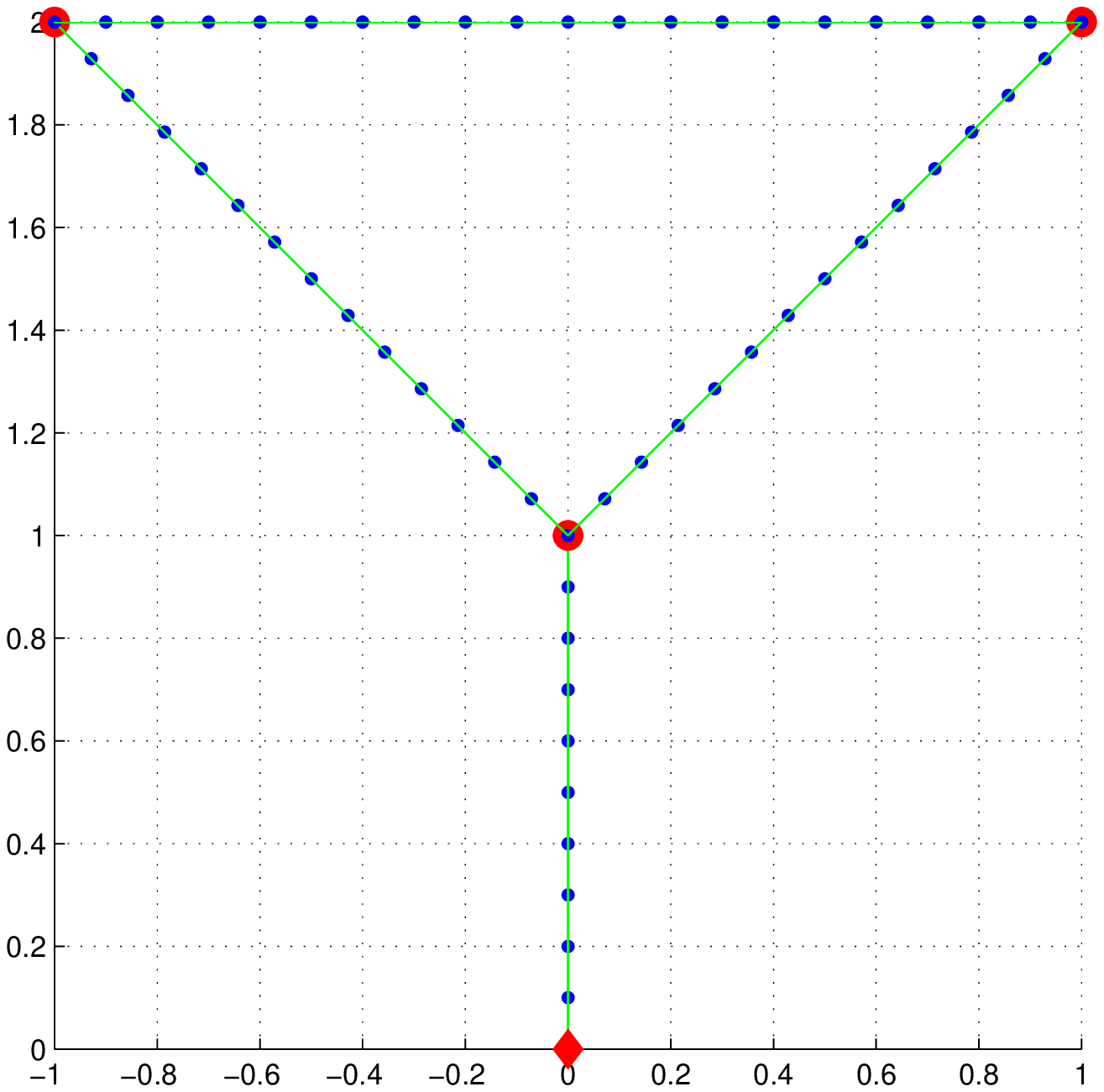,width=4.5cm}\hspace{2cm}\epsfig{figure=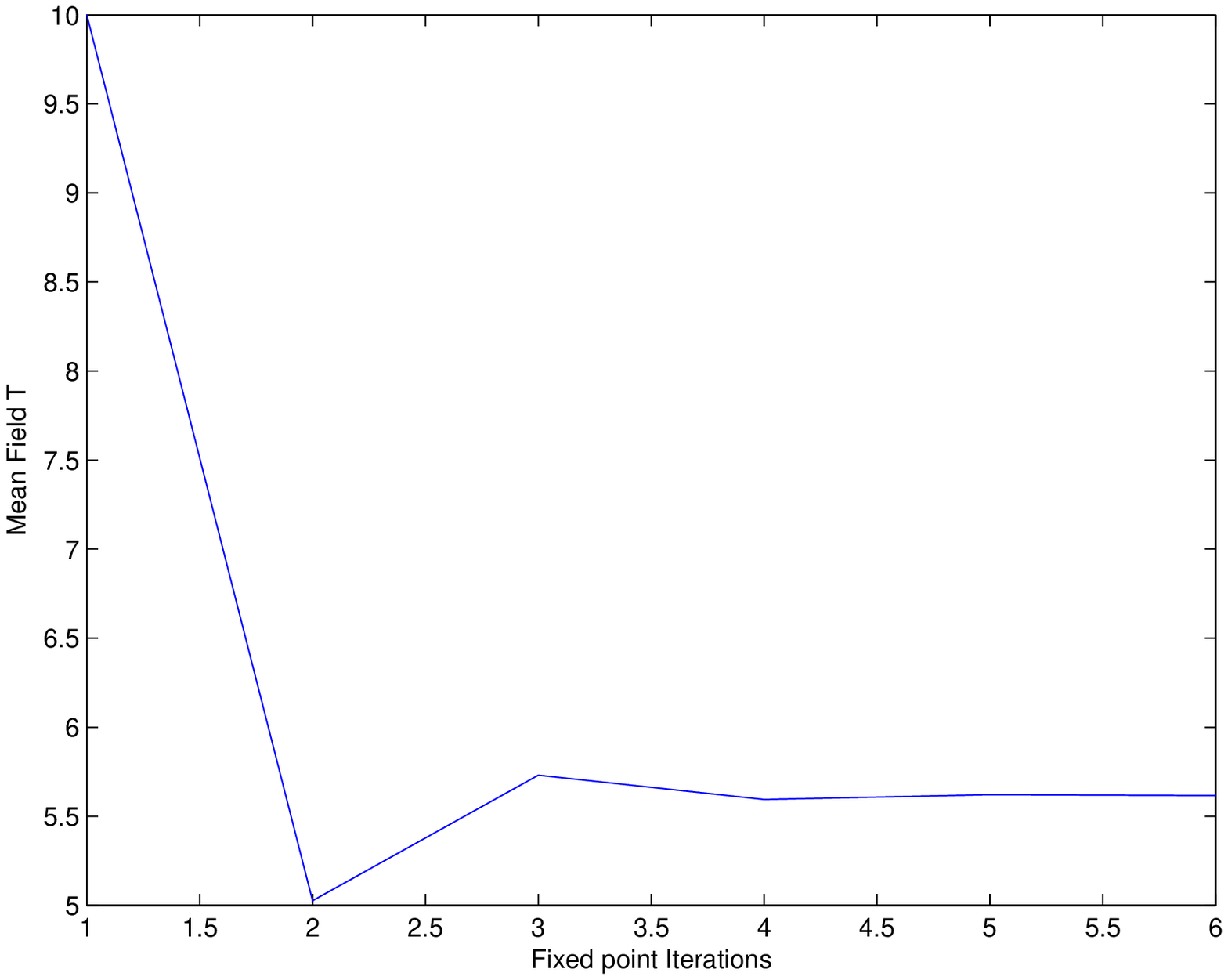,width=6.5cm}
\caption{ {{ Left: Graph configuration. Right: mean field approximated time $T_2$ vs. fixed point iterations, computed by Algorithm \ref{algo} }}  }
\label{Test1grafo}
\end{center}
\end{figure}

\begin{figure}[ht!]
\begin{center}
\epsfig{figure=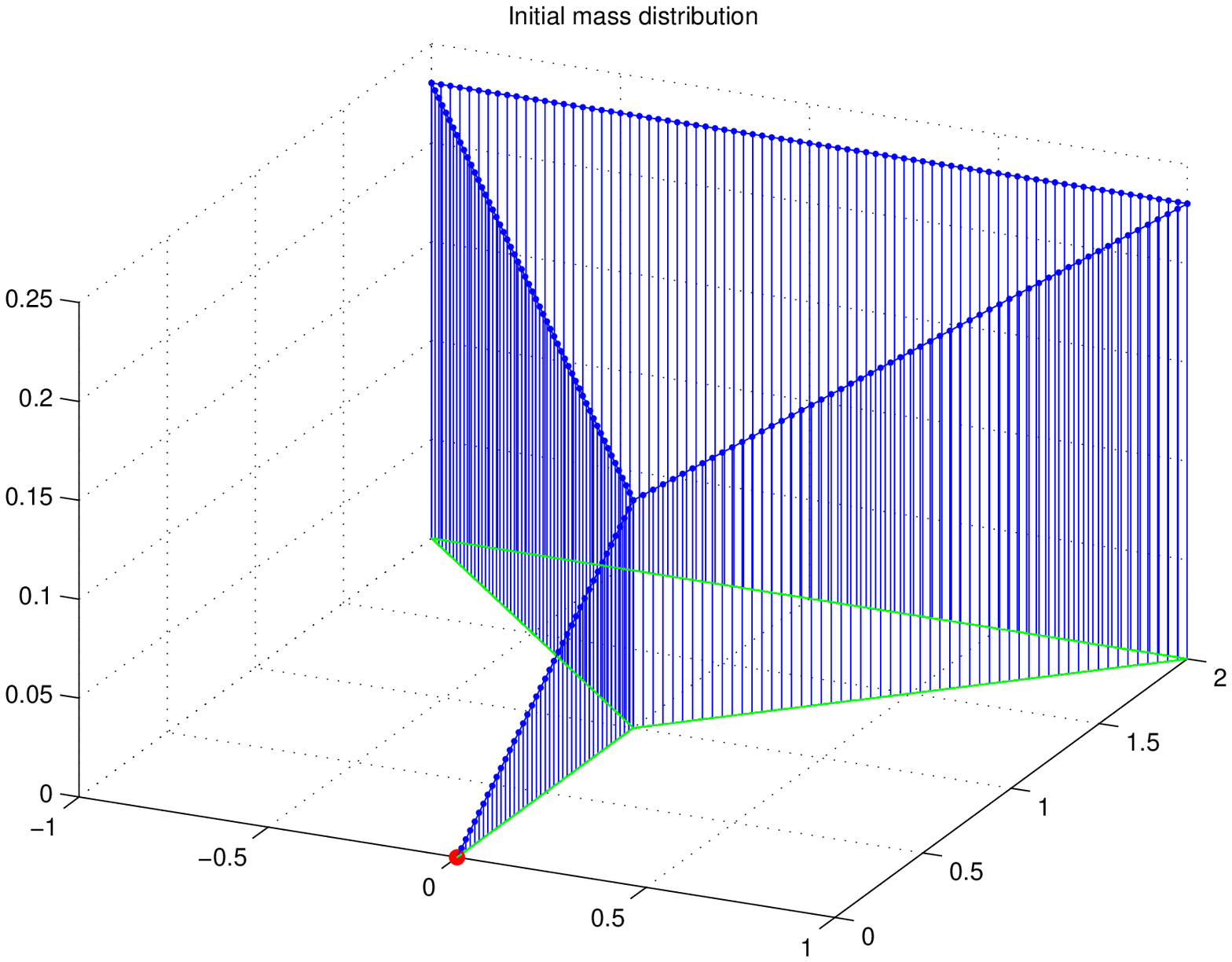,width=5cm}\epsfig{figure=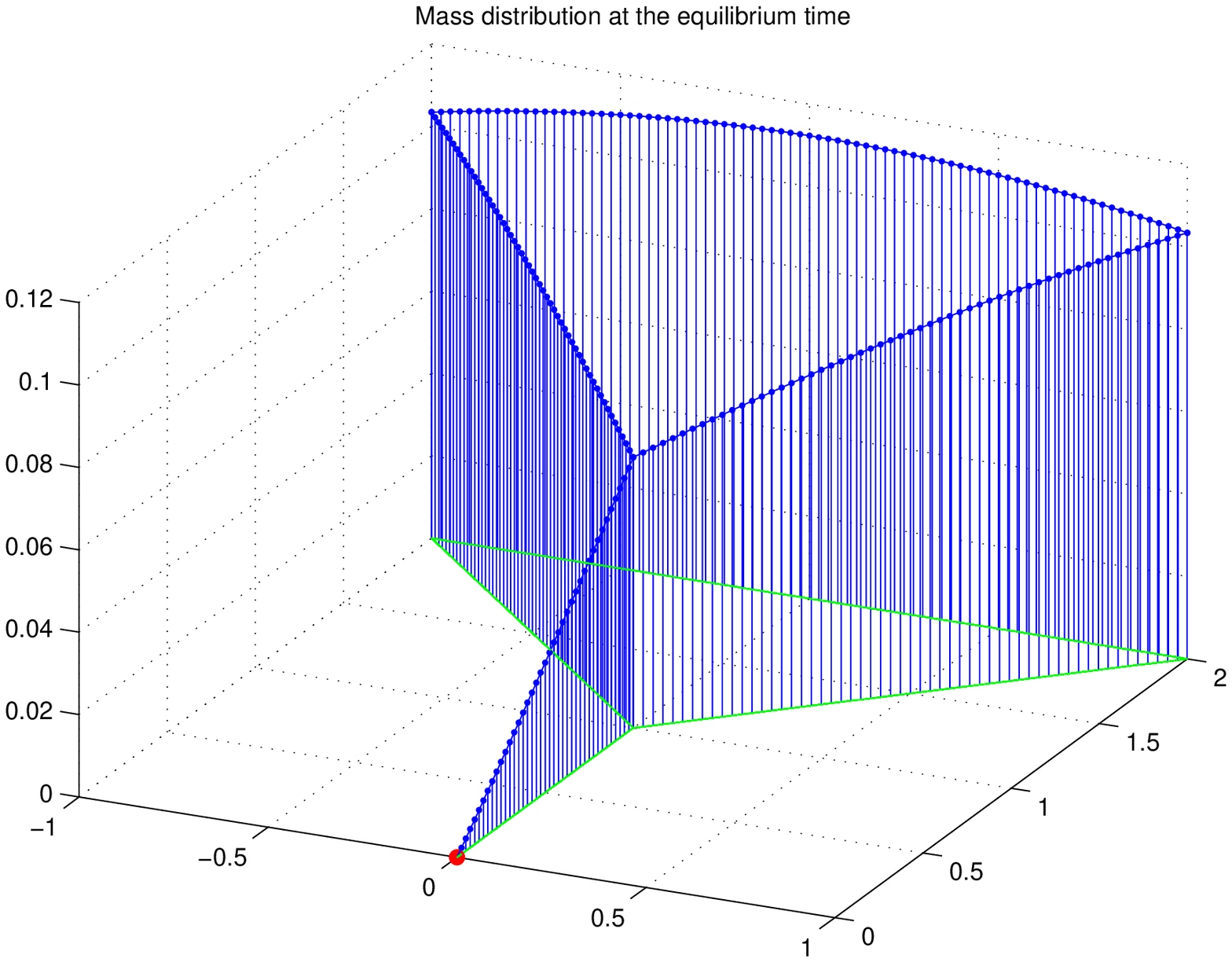,width=5cm}\epsfig{figure=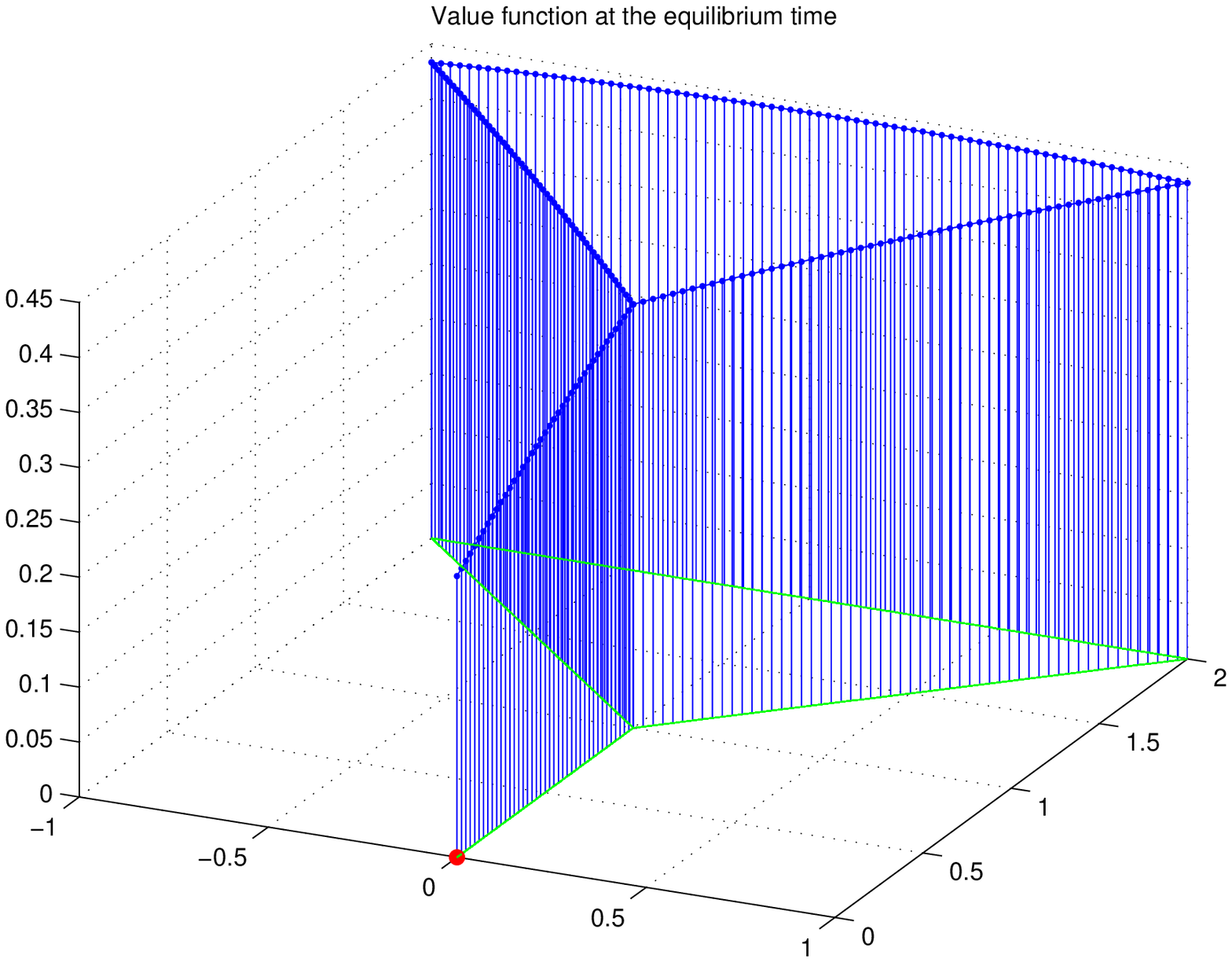,width=5cm}
\caption{ {{ Test 1: Initial mass distribution (left), distribution at the equilibrium time (center), value function (right) }}  }
\label{Test1mass}
\end{center}
\end{figure}

We consider a simple graph with four vertexes and four edges, as shown in Fig.\ref{Test1grafo}.
The initial mass distribution is given by
$$m_0(x)=\frac{g(x)}{\int_{\Gamma}g(y) d y},$$
where $g(x)$ is the restriction to $\G\subset \R^2$ of the function $|x|$.
The scheduled time is $t_0=0.5$, the maximal waiting time is $T_{max}=10$, %the viscosity parameter is $\nu=1$,
the  cost function is $$c_T(s)=0.1 \max(s-t_0,0)+0.1 \max(T-s,0)$$
and the percentage value of the  expected players is $\theta=0.5$.\\
For each arc $j\in J$, we consider  the same space step $h_j=h$ and we run a series of numerical tests  varying the space step according to the first   column of Table \ref{tab:test1}. The time step has to verify the stability condition $\Delta t <h^2$ and then we  choose $\Delta t=h^2/4$.
For each test we compute the following error
\begin{equation}\label{err}
E_h(T)=\left|1-\theta-\sum_j\sum_i \psi^N_{j,i}\phi^N_{j,i} h_j \right|\simeq \left|1-\theta-\int_{\Gamma}m(x,T)dx\right|,
\end{equation}
where $N$ is such that $T_2=N\Delta t$ in Algorithm \ref{algo}.
%which it measure the how much mass has been lost in the numerical computation.
Since $\theta$ represents the percentage of player exited from the boundary vertex $v_0$,  then $1-\theta$  represents the percentage of  the residual  population and the term on the right side  of \eqref {err} should be zero. This error is shown in the second column  of Table \ref{tab:test1}.  In the third and fourth columns we show the computed mean time $T_2$, and the  number of iterations needed by the Algorithm \ref{algo} to converge when $\tau=10^{-4}$ and $T_1=10$.
Table  \ref{tab:test1} shows small values for $E_h(T)$ and, even if we do not observe a monotone behavior, the smallest value is attained with the finer space grid. \\
The graph on the right of Figure \ref{Test1grafo}  shows the convergence of the approximated mean field time $T_2$, computed by Algorithm \ref{algo} with space step $h=2.50  \cdot 10^{-2}$. On the horizontal axis are the iterations of the fixed point, while on
vertical axis the corresponding approximated mean field time $T_2$.
In Fig.\ref{Test1mass}, we show  the initial mass distribution (left), equilibrium mass distribution (center)  and the corresponding value function (right).
\begin{table}[ht!]
\begin{center}
\label{tab:test1}
\begin{tabular}{|c|c|c|c|}\hline
$h$ & $E_h(T_2)$ & $T_2$& iterations \\
\hline \hline
$1.00 \cdot 10^{-1}$& $8.27 \cdot 10^{-4} $&5.687& 6\\ \hline
$5.00 \cdot 10^{-2}$ & $1.34 \cdot 10^{-3}$ & $ 5.639$&  7 \\ \hline
$2.50  \cdot 10^{-2} $& $9.04 \cdot 10^{-4}$ & $ 5.617$&  8 \\ \hline
$1.25  \cdot 10^{-2}$& $5.02 \cdot 10^{-4}$ & $ 5.622$&  6 \\ \hline
\end{tabular}
\end{center}
\caption{Space steps (first column), $E_h(T_2)$ defined in \eqref{err} (second column), approximated  mean field $T_2$ (third column), number of fixed point iterations (last column)}
\end{table}

\subsection{Example 2: A more general graph}
\begin{figure}[ht!!]
\begin{center}
\epsfig{figure=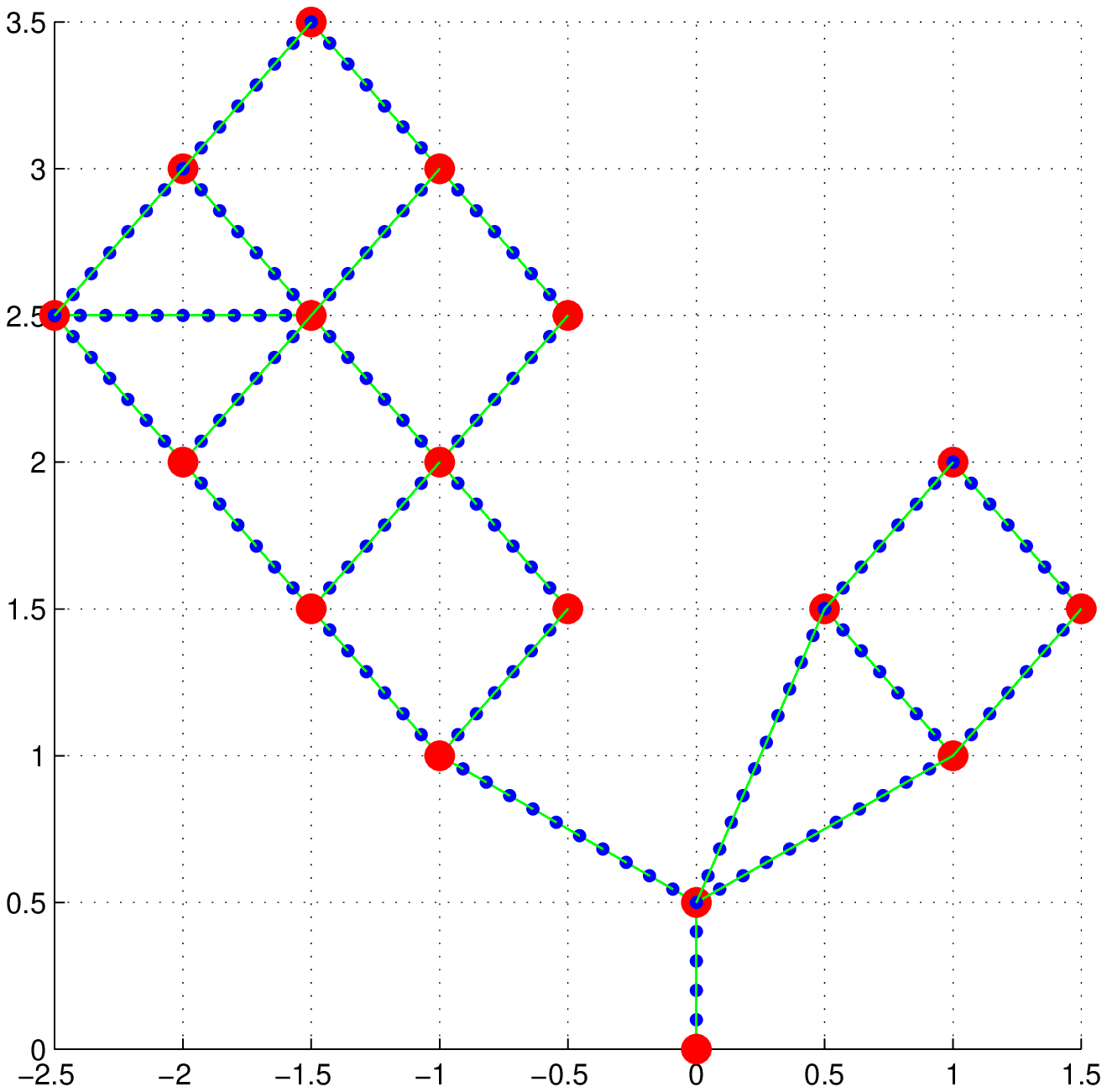,width=4.5cm}\hspace{2cm}\epsfig{figure=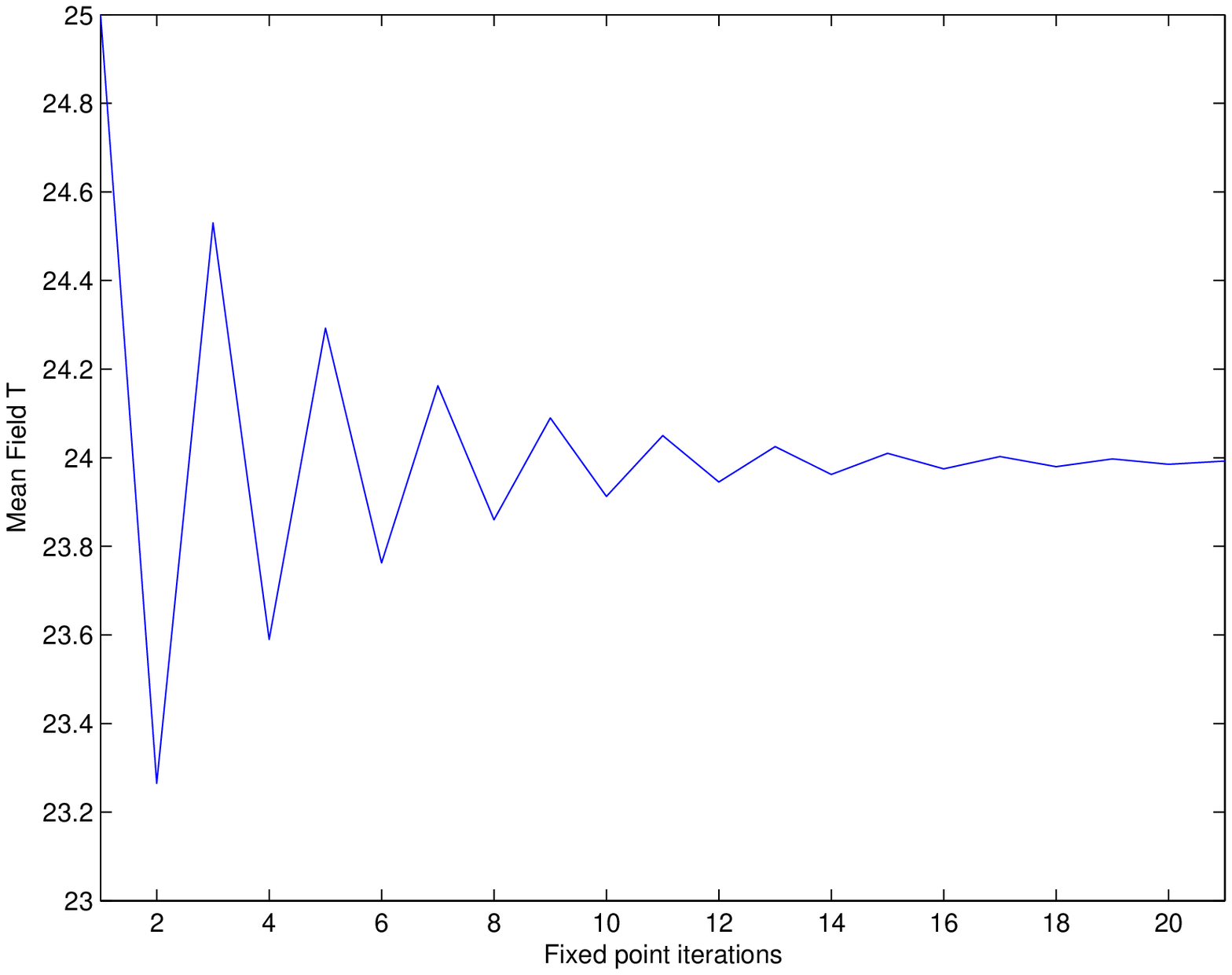,width=6.5cm}
\caption{ {{ Left: Graph configuration. Right: approximated mean field  time $T_2$ vs. fixed point iterations, computed by Algorithm \ref{algo} }}  }
\label{Test2grafo}
\end{center}
\end{figure}

\begin{figure}[ht!]
\begin{center}
\epsfig{figure=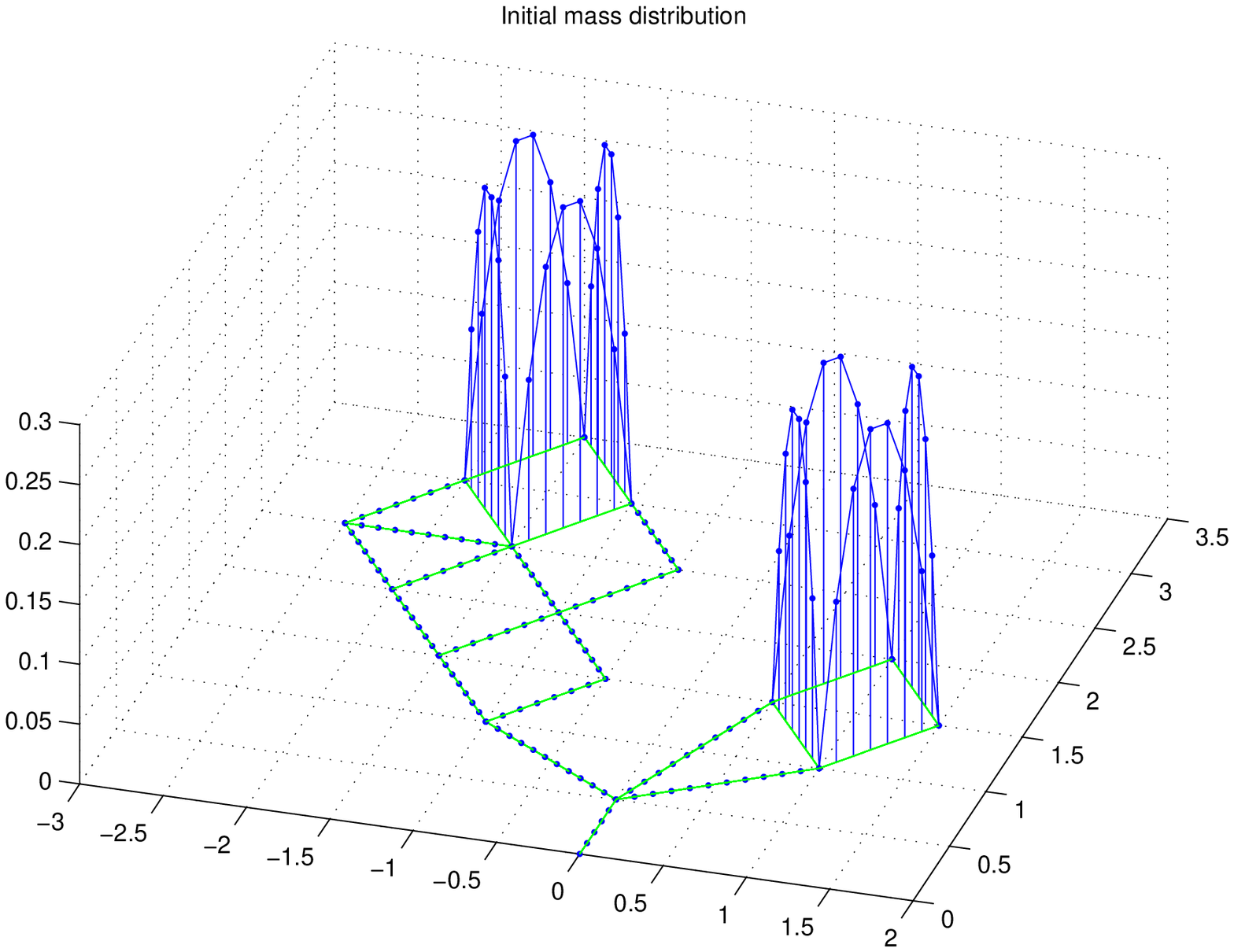,width=5cm}\epsfig{figure=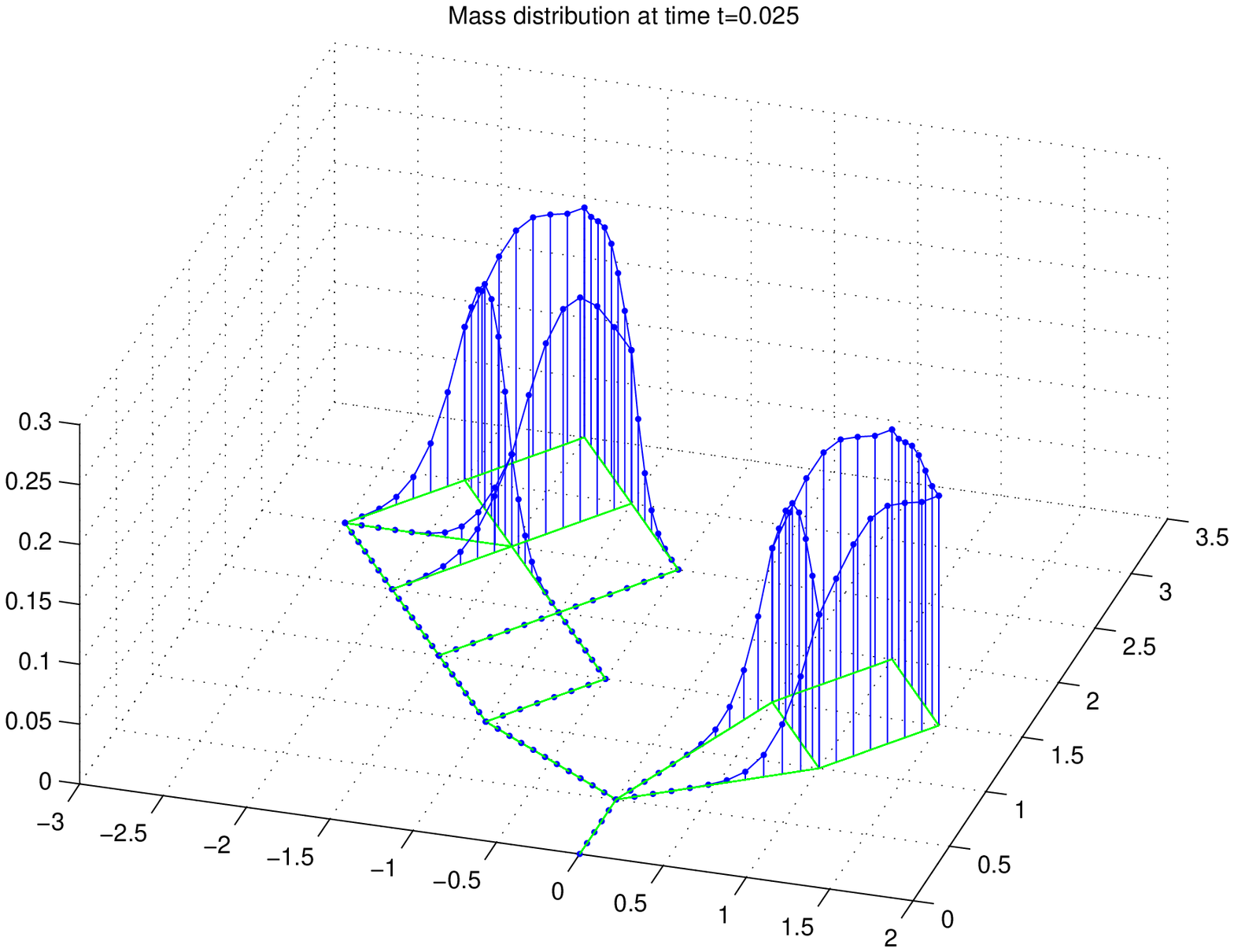,width=5cm}\epsfig{figure=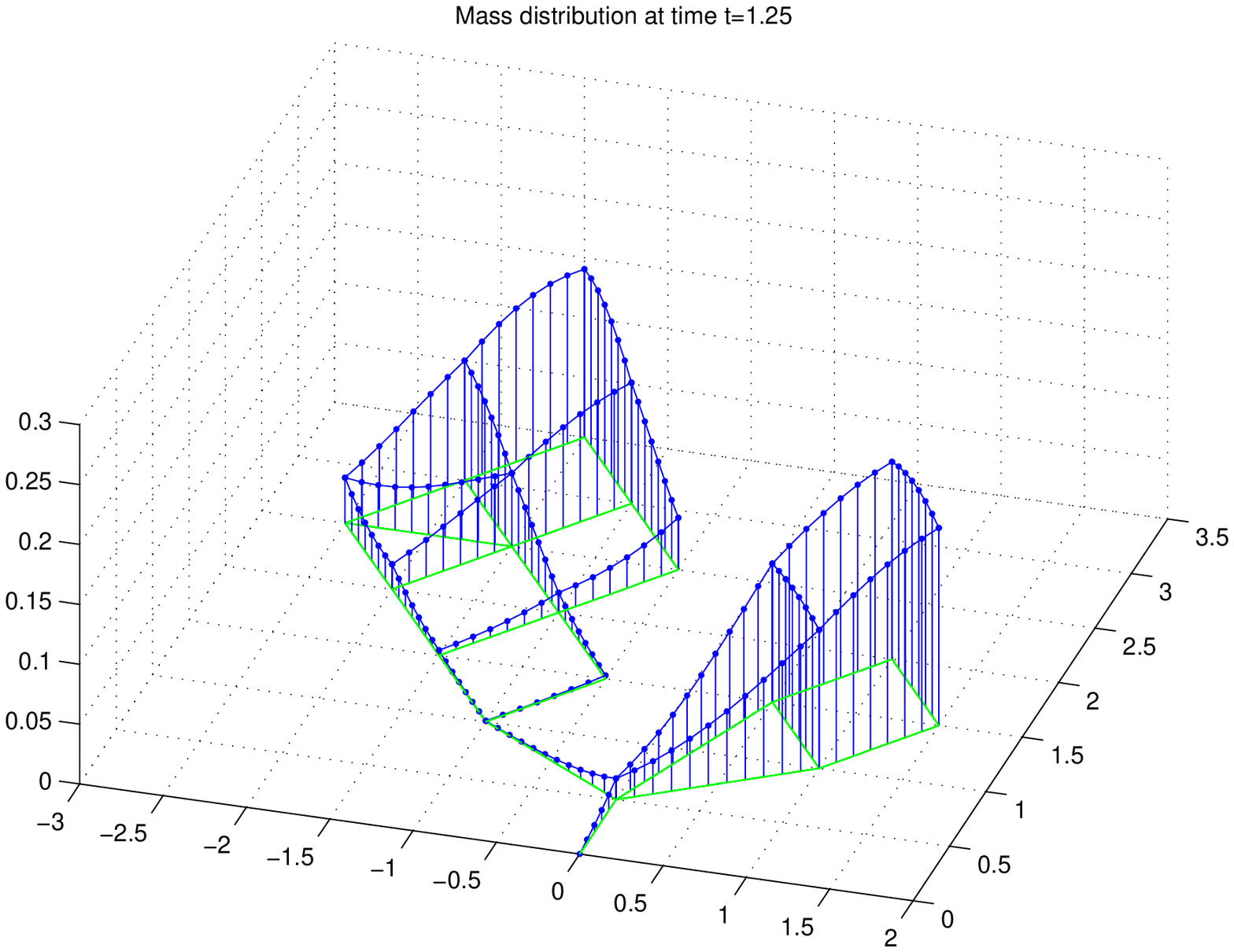,width=5cm}\\
\epsfig{figure=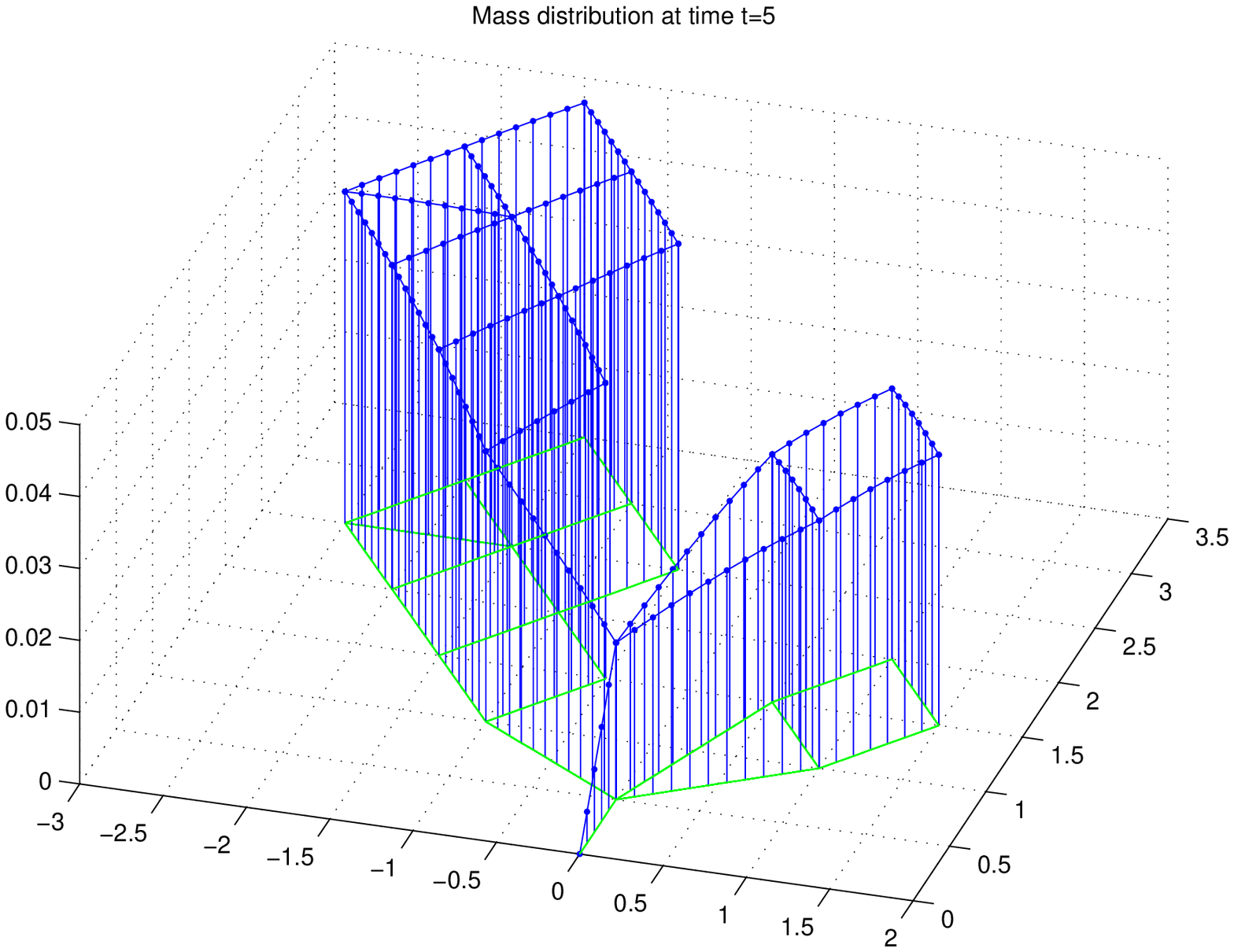,width=5cm}\epsfig{figure=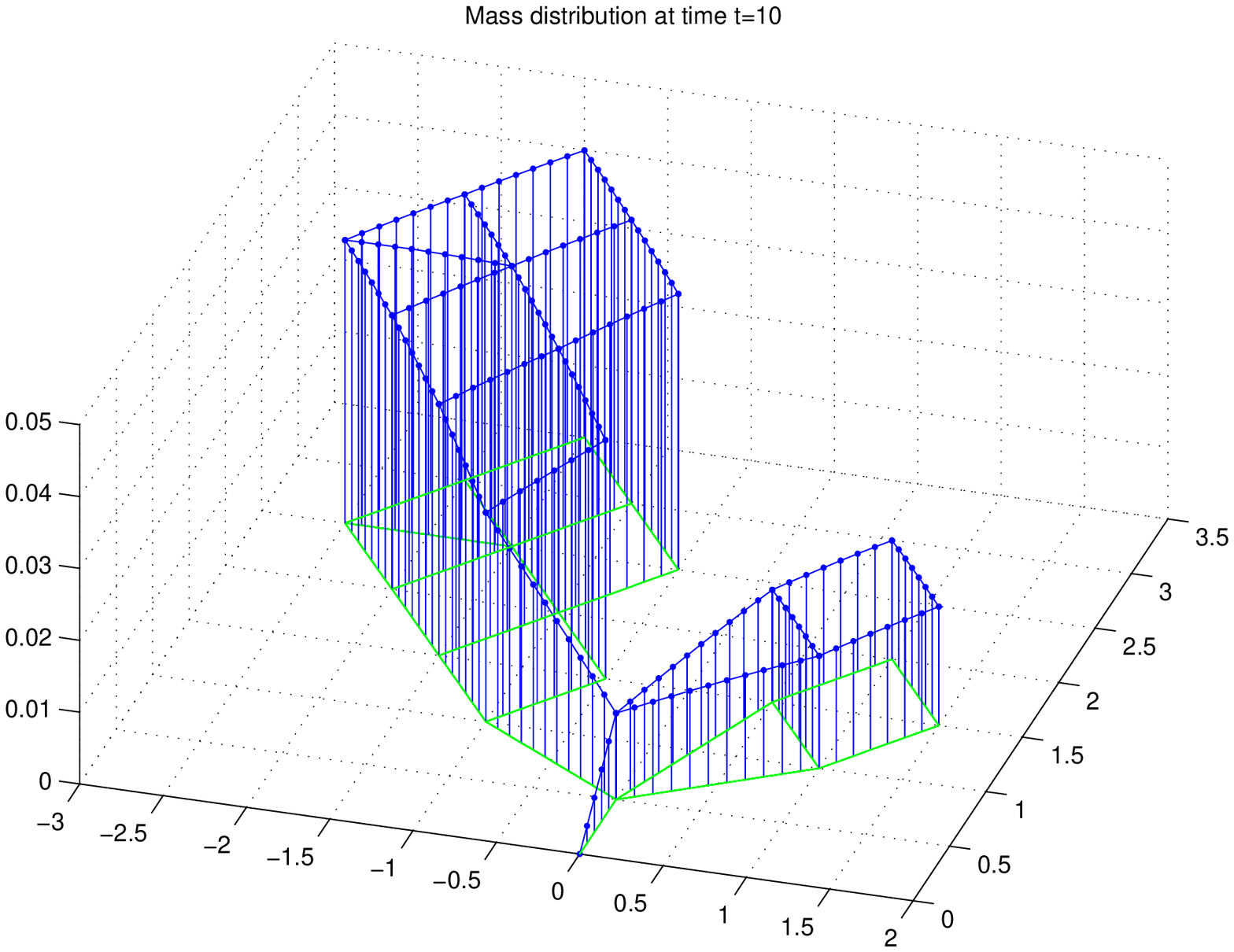,width=5cm}\epsfig{figure=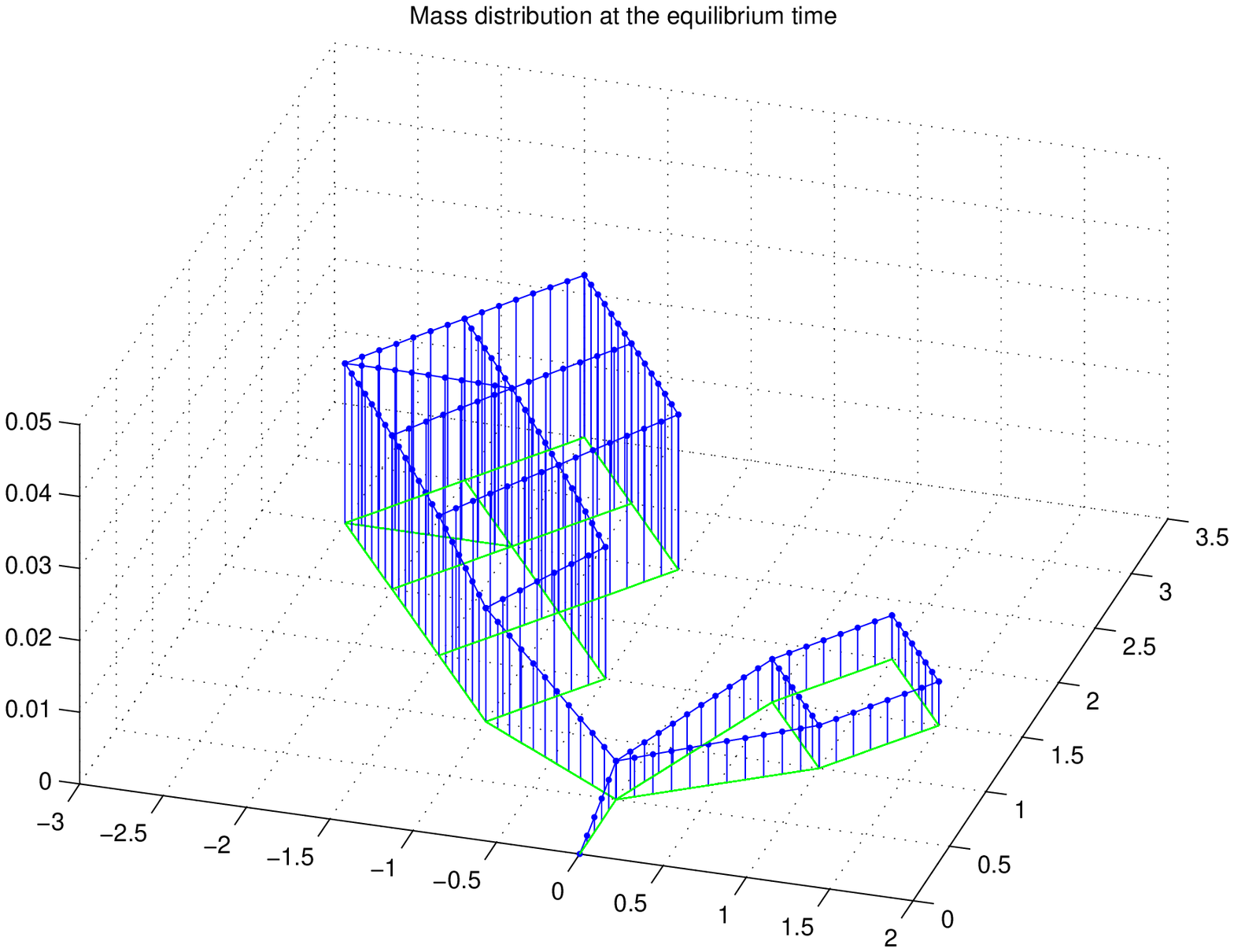,width=5cm}
\caption{ Test 2: Mass distribution at time: $t=0, 0.025,1.25, 5, 10, T=24$}
\label{Test2mass}
\end{center}
\end{figure}

We consider a more general graph with 17 vertexes and 22 edges, see   Fig.\ref{Test2grafo}.
The initial mass distribution is given by
$$m_0(x)=\frac{g(x)}{\int_{\Gamma}g(y) d y}, \quad g(x)=\max(0.5-|x-p_1|^2,0)+\max(0.5-|x-p_2|^2,0)\quad x\in \Gamma,$$
with $p_1=(1,3/2)$ and $p_2=(-3/2,3)$.  It describes the distribution of two populations, one   concentrated around the point $p_1$, the other one around $p_2$.\\
The scheduled time is $t_0=0.5$, the maximum waiting time is $T_{max}=25$, %the viscosity parameter is $\nu=1$,
the  cost function
$$c(s)=0.1 \max(s-t_0,0)+0.1 \max(T-s,0)$$
and the  expected percentage of arrival players is $\theta=0.7$. The Algorithm \ref{algo} is run with $h=0.05$, $\Delta t=\frac{h^2}{4}$ and $\tau=0.05$. We get $T=23.99$ with error  $E_h(T)=2.35 \cdot 10^{-2}$.
The graph on the right of  Figure \ref{Test2grafo} shows the convergence of the approximated mean field time $T_2$ computed by Algorithm \ref{algo}: on the horizontal axis is the    number of iterations of the fixed point algorithm, whereas on
the vertical axis the corresponding mean field time.\\
%\begin{figure}[ht!]
%\begin{center}
%\epsfig{figure=figure/Tmeanconv.eps,width=5cm}
%\label{Test2Tmean}
%\end{center}
%\end{figure}
Figure \ref{Test2mass} shows the mass evolution at different times.
It can be observed that at the initial time the diffusion spreads   the population in all the   directions on the graph,
later the cost \eqref{costtoy} favors the population  closer to $v_0$  to reach the exit before of the population farther away.

%%%%%%%%%%%%%%%%%
%               %
%%%%%%%%%%%%%%%%%
\section{Appendix}
\label{Sec:Appendix}
\begin{Proofc}{Proof of Prop. \ref{exisFP0}}
For the sake of simplicity, $K$ will denote a constant independent of $\mu_0$ and $w$ and it may change from line to line.
Invoking \cite[Theorem 5.4]{mu} (see also: \cite[Theorem 3.2]{fw}, \cite[Theorem 3.6]{kms} or \cite[Theorem 5.8]{mr}) we obtain that there exists a unique classical solution $\mu$ to problem \eqref{FPnet} which fulfills the estimate
\begin{equation}\label{FPMugn}
|\mu|_\infty\leq K |\mu_0/w(\cdot,0)|_\infty.
\end{equation}

For $\tilde e:=\pi_0([l_0/4,3l_0/4])$, we claim that $\mu$ belongs to $C^{(2+\a,1+\a/2)}(\tilde e\times [0,T_{max}])$ with
\begin{equation}\label{33lip}
|\mu|^{(2+\a,1+\a/2)}_{\tilde e \times [0,T_{max}]}\leq K|\mu_0/w(\cdot,0)|^{(2+\a)}_{\G}.
\end{equation}
In order to prove this estimate, we introduce two families of functions $\{\tilde \mu_{0,n}\}_n$ and $\{\tilde \mu_{1,n}\}_n$ such that
\begin{equation*}
\begin{array}{c}
\tilde \mu_{0,n}, \tilde \mu_{1,n}\in C^1([0,T_{max}]),\qquad
|\tilde \mu_{0,n}|_\infty+|\tilde \mu_{1,n}-\mu(v_1,\cdot)|_\infty\to 0\quad \textrm{as }n\to+\infty,\\
\tilde \mu_{0,n}(0)=0,\qquad \tilde \mu_{0,n}'(0)=D^2\left(\frac{\mu_0(\cdot)}{w(\cdot,0)}\right)(v_0),\\
\tilde \mu_{1,n}(0)=\frac{\mu_0(v_1)}{w(v_1,0)},\qquad \tilde \mu_{1,n}'(0)= D^2\left(\frac{\mu_0(\cdot)}{w(\cdot,0)}\right)(v_1),\\
\end{array}
\end{equation*}
By standard regularity theory for parabolic equations on domains in Euclidean spaces, the problem
\begin{equation*}
\ds\left\{
\begin{array}{ll}
\pd_t \mu_n-\pd^2_x\mu_n=0\quad& (x,s)\in (0,l_0)\times (0,T_{max})\\[6pt]
\mu_n(0,s)=\tilde\mu_{0,n}(s),\quad \mu_n(l_0,s)=\tilde \mu_{1,n}(s) & s\in[0, T_{max}] \\[6pt]
\mu_n(x,0)= \frac{\mu_0(x)}{w(x,0)} & x\in(0,l_0)
\end{array}\right.
\end{equation*}
admits a unique classical solution $\mu_n$ which belongs to $C^{(2+\a, 1+\a/2)}((0,l_0)\times (0,T_{max}))$ for some $\a$ depending only on the features of the equation. By \cite[Theorem IV.10.1]{lsu}, we deduce the following estimate in the domain $(l_0/4,3l_0/4)\times (0,T_{max})$
\begin{equation*}
|\mu_n|^{(2+\a,1+\a/2)}_{(l_0/4,3l_0/4)\times (0,T_{max})}\leq K\left(|\mu_0/w(\cdot,0)|^{(2+\a)}_{e_0} + |\mu_n|_\infty \right).
\end{equation*}
By Ascoli theorem, as $n\to +\infty$, (eventually, passing to a subsequence), the function $\mu_n$ converges uniformly to some function $v$ and the same happens for $\pd_t \mu_n$, $\pd_x \mu_n$ and $\pd_{x}^2 \mu_n$ with the corresponding derivatives of $v$. By the stability result we get $v=\mu$. Moreover, passing to the limit in the last estimate, we obtain
\begin{equation*}
|\mu|^{(2+\a,1+\a/2)}_{(l_0/4,3l_0/4)\times (0,T_{max})}\leq K\left(|\mu_0/w(\cdot,0)|^{(2+\a)}_{e_0} + |\mu|_\infty \right)
\end{equation*}
and, taking into account estimate \eqref{FPMugn} and the definition of the sub-edge $\tilde e$, we accomplish the proof of claim \eqref{33lip}.

We observe that the function $\mu_{|e_{0,1/2}\times(0,T_{max})}$ is the unique classical solution to problem
\begin{equation*}
\ds\left\{
\begin{array}{ll}
\pd_t \bar \mu-\pd^2_x\bar \mu=0& (x,s)\in e_{0,1/2}\times (0,T_{max})\\[6pt]
\bar \mu(v_0,s)=0,\quad \bar \mu(v'_{1/2},s)=\mu(v'_{1/2},s)& s\in[0, T_{max}] \\[6pt]
\bar \mu(x,0)= \frac{\mu_0(x)}{w(x,0)} & x\in e_{0,1/2}
\end{array}\right.
\end{equation*}
which is a standard initial-boundary value problem on an Euclidean domain.
Invoking \cite[Theorem IV.9.1]{lsu}, we infer that, for every $q\geq 1$, $\mu$ belongs to $W^{2,1}_{q, e_{0,1/2}\times [0,T_{max}]}$ with
\begin{equation*}
|\mu|^{2,1}_{q, e_{0,1/2}\times [0,T_{max}]}\leq K
\left(|\mu_0/w(\cdot,0)|^{(2)}_{e_{0,1/2}}+|\mu(v'_{1/2},\cdot)|^{(1)}_{(0,T_{max})}\right).
\end{equation*}
Owing to \eqref{33lip}, estimate \eqref{33weak} is achieved.
\end{Proofc}
%\begin{Remark}
%Let us observe that, in the previous proof, estimate \eqref{33lip} is used only for the Lipschitz continuity of $\mu(v'_{1/2},\cdot)$.
%However, the proof of this regularity is almost straightforward provided that $\mu_0$ filfills the Kirchhoff condition:
%\begin{equation*}
%\sum_{j\in Inc_v} a_{ij}\rho_{ij} \pd_{x,j} \mu_0(v_i ,s)=0\qquad (v_i,s)\in \G_T\times (0,T_{max}).
%\end{equation*}
%Actually, in this case, the functions $v^\pm(x,t):=\pm |\pd_x^2\mu_0/w(\cdot,0)|_\infty t + \mu_0(x)/w(x,0)$ are respectively a super- and a subsolution to problem \eqref{FPnet}; it is worth to observe that they verify the Kirchhoff condition because of previous relation and Kirchhoff condition for $\phi$. Whence, the comparison principle ensures
%\begin{equation*}%\label{33weak2}
%|\mu(x,t)|\leq |\pd_x^2\mu_0/w(\cdot,0)|_\infty t + \mu_0(x)/w(x,0)\qquad \forall (x,t)\in \bar \G\times [0,T_{max}].
%\end{equation*}
%We introduce the functions $\mu^\pm_h(x,t):=\mu(x,t+h)\pm |\pd_x^2\mu_0/w(\cdot,0)|_\infty h$.
%Taking advantage of the previous inequality, these functions are respectively a super- and a subsolution to problem \eqref{FPnet}.
%Invoking again the comparison principle, we deduce
%\[|\mu(x,t)-\mu(x,t+h)|/h\leq |\pd_x^2\mu_0/w(\cdot,0)|_\infty. \]
%\end{Remark}

\begin{Proofc}{Proof of Prop. \ref{exisFP}}
We shall improve some arguments of the proof of Proposition \ref{exisFP0} taking advantage of the stronger compatibility condition given by \eqref{cmpv0}. Here, the constant $K$ is independent of $\mu_0$ and $w$ and it may change from line to line.

%QUESTA PROPRIETA' DI LIPSCHITZIANITA' FORSE E' INUTILE: We claim that the function $\mu(v'_{3/4}, \cdot)$ is Lipschitz continuous on $[0,T_{max}]$ with a Lipschitz constant less or equal to $e^{T_{max}}\|\mu\|_\infty$. Actually, in order to prove this property, it suffices to observe that the functions
%\begin{equation*}\mu^\pm_h:=\mu(\cdot,\cdot+h) \pm 2h\|\mu_0/w(\cdot,0)\|^{(2)}_\G\end{equation*}
%are respectively a super- and a subsolution to the staring problem. FINE PROPRIETA' DI LIPSCHITZIANITA'.

We consider the family of functions $\{\tilde\mu_{1,n}\}_n$ introduced in the proof of Proposition \ref{exisFP0}.
%We consider a family of functions $\{\tilde\mu_n\}_n$ such that
%\begin{equation*}\begin{array}{ll}
%\tilde\mu_n\in C^1([0,T_{max}])&\quad\|\tilde\mu_n -\mu(v'_{3/4},\cdot)\|_\infty\to 0\quad \text{as }n\to+\infty\\
%\tilde\mu_n (0)=\frac{\mu_0(v'_{3/4})}{w(v'_{3/4},0)},&\quad\tilde\mu_n' (0)=\nu D^2\left(\frac{\mu_0(\cdot)}{w(\cdot,0)}\right)(v'_{3/4}).
%\end {array}\end{equation*}
By standard regularity theory for parabolic equations on domains in Euclidean spaces, the problem
\begin{equation*}
\ds\left\{
\begin{array}{ll}
\pd_t \mu_n-\pd^2_x\mu_n=0\quad& (x,s)\in (0,l_0)\times (0,T_{max})\\[6pt]
\mu_n(0,s)=0,\quad \mu_n(l_0,s)=\tilde \mu_{1,n}(s)& s\in[0, T_{max}] \\[6pt]
\mu_n(x,0)= \frac{\mu_0(x)}{w(x,0)} & x\in(0,l_0)
\end{array}\right.
\end{equation*}
admits a unique classical solution $\mu_n$ which belongs to $C^{(2+\a, 1+\a/2)}((0,l_0)\times (0,T_{max}))$ for some $\a$ depending only on the features of the equation. By \cite[Theorem IV.10.1]{lsu}, we deduce the following estimate in the domain $(0,l_0/2)\times (0,T_{max})$
\begin{equation}\label{iv101}
|\mu_n|^{(2+\a,1+\a/2)}_{[0,l_0/2]\times [0,T_{max}]}\leq K\left(| \mu_0/w(\cdot,0)|^{(2+\a)}_{[0,l_0]} + |\mu_n|_\infty \right).
\end{equation}
By Ascoli theorem, as $n\to +\infty$, (eventually, passing to a subsequence), the function $\mu_n$ converges to some function $v$ uniformly and the same happens for $\pd_t \mu_n$, $\pd_x \mu_n$ and $\pd_{x}^2 \mu_n$ with the corresponding derivatives of $v$. By the stability result we get $v=\mu$. Moreover, passing to the limit in the estimate \eqref{iv101}, we obtain
\begin{equation*}
|\mu|^{(2+\a,1+\a/2)}_{[0,l_0/2]\times [0,T_{max}]}\leq K\left(| \mu_0/w(\cdot,0)|^{(2+\a)}_{[0,l_0]} + |\mu|_\infty \right).
\end{equation*}
Finally, taking into account estimate \eqref{FPMugn}, we accomplish the proof.\par
The second part of the  statement is a consequence of \cite{vb}; actually, in this case, the compatibility conditions are ensured by \eqref{FPallaVB}.
Invoking  \cite{vb}, we obtain
\begin{equation*}
|\mu|^{(2+\a,1+\a/2)}_{\G\times[0,T_{max}]}\leq K_0 |\mu_0/w(\cdot,0)|^{(2+\a)}_{\G}
\end{equation*}
where $K_0$ is the same constant as in Proposition \ref{exisheat}.
\end{Proofc}
\begin{Remark}
As one can easily check, in the proof of previous Proposition \ref{exisFP}, hypothesis \eqref{cmpv0} is needed only for guaranteeing the compatibility condition in $v_0$. As a matter of fact, it can be replaced by: $\pd_x^2(\mu_0(\cdot)/w(\cdot,0))(v_0)=0$.
\end{Remark}

 %%%%%%%%%%%%%%%%%%%%%%%%%%%%%%%%%%%%%%%%%%%%%%
\textbf{Acknowledgment.}The authors wish to thank Adriano Festa for his help in the algorithm for programming the graph structure

\end{document}